\documentclass{article}
\usepackage{arxiv}

\usepackage[utf8]{inputenc} 
\usepackage[T1]{fontenc}    
\usepackage{url}            
\usepackage{booktabs}       
\usepackage{nicefrac}       
\usepackage{microtype}      

\usepackage[english]{babel}
\usepackage{caption}
\usepackage{subcaption}
\usepackage{graphicx}

\usepackage{amsfonts,amssymb,amsmath,amsthm}
\usepackage{breqn}
\usepackage{epsfig}
\usepackage{epstopdf}
\usepackage{amstext}
\usepackage{amsbsy}
\usepackage{amsopn}
\usepackage{float}
\usepackage{lipsum}
\usepackage{appendix}
\usepackage{multirow}
\usepackage{hyperref}
\usepackage{cleveref}
\usepackage{setspace}

\usepackage{url}
\usepackage{color}

\newcommand{\norm}[1]{\left\lVert#1\right\rVert}

\newcommand\pd[3]{\frac{\partial^{#1} #2}{\partial #3^{#1}}}
\newcommand\status{false}

\title{Natural response of non-smooth oscillators using homotopy analysis combined with Galerkin projections}

\author{
  Jeet Desai\\
  Universit\'e Pierre et Marie Curie\\
  Paris, France \\
  \texttt{jeet.desai@etu.upmc.fr } \\
   \And
 Amol Marathe \\
  Department of Mechanical Engineering\\
  BITS-Pilani\\
  Pilani, Rajasthan \\
  \texttt{amolmm@pilani.bits-pilani.ac.in } \\}

\begin{document}

\maketitle

\begin{abstract}
We propose homotopy analysis method in combination with Galerkin projections to approximate the natural response of non-smooth oscillators with discontinuities of type Heaviside, signum, modulus etc. While constructing the homotopy, we think of convergence-control parameter as a function of the embedding parameter and call it a convergence-control function. Homotopy analysis provides an expression for the natural frequency of the oscillator that also includes free parameters arising from the convergence-control function. Generating extra equations using Galerkin projections and solving the same numerically gives the approximate natural response of the non-smooth oscillators. We also seek aperiodic natural response of a unilaterally constrained simple pendulum. The superiority of the approach is well-established over the method of harmonic balance, Lindstedt-Poincaré method and conventional homotopy analysis method.
\end{abstract}


\section{Introduction}
\label{intro}

Oscillatory motion when constrained rigidly is typically governed by a second order non-smooth differential equation. Several problems from mechanical engineering, e.g., vibrations of a spring-mass system with unequal restraints, pendulum with impact, a gear-pair with backlash and friction etc. are modelled using second order differential equations involving discontinuous mathematical functions such as signum, Heaviside, modulus etc. Natural response of such oscillators may be obtained with ease using numerical integration. Also, several perturbation-like methods such as parameter expansion, homotopy perturbation, modified Lindstedt-Poincar\'{e}, variational iteration have been applied successfully to get the periodic solution as well as the approximate analytical estimate of the natural frequency \cite{HE2002309,XU2007148}. In the context of non-smooth oscillators, such methods are typically applied to find periodic solutions \cite{LIU2005577,WANG2008688}. For example, the homotopy perturbation method has been successfully applied to obtain the solution of a nonlinear oscillator where the discontinuity is modelled using the modulus function \cite{HE2004287}. Variational iteration technique (VIM) is applied to several strongly nonlinear oscillators including the ones involving discontinuities \cite{HE1999699,RAFEI2007614,SHOU20092416}. Frequency of the limit cycles of several non-smooth oscillators is determined using an artificial parameter-Linstedt–Poincaré method \cite{RAMOS2008738}. Non-periodic response where amplitude decays with time or frequency does not remain constant are harder to approximate using such methods. One such example is a unilaterally constrained simple pendulum where rigid body collision is modelled using Newton's impact law. The constraint is usually not accommodated in the governing equation, thereby restricting the applicability of such methods.

Here, we develop a framework that combines homotopy analysis method (HAM) with Galerkin projections to obtain not only the periodic response of certain non-smooth oscillators, but also that of the unilaterally constrained pendulum where the solutions are not periodic. The deformation is parametrized using an embedding parameter $p\in[0,1],$ deforming the solution of a linear operator ($p=0$) till it reaches that of the non-smooth oscillator or the desired solution ($p=1$). To approximate the periodic responses, the choice of the linear operator is a simple harmonic oscillator. While obtaining the periodic response via HAM, we construct the solution at every order over only first period and also remove the secular terms over the same. To accelerate the convergence of the series solution at $p=1,$ we introduce a convergence-control function $h(p),$ conventionally considered to be a constant. By following the homotopy procedure, different derivatives of this function evaluated at $p=0$ appear naturally as unknowns in the expressions for the periodic solution as well as the natural frequency. When applied up to $n$-th order deformation equation, our method produces $n$ unknowns compared to HAM which produces just one ($h$). To determine these, we apply Galerkin projections with appropriate weighting functions and generate extra equations, usually non-algebraic in the context of non-smooth oscillators. Solving the system of equations numerically yield the desired periodic response and the natural frequency. However, this last step makes the procedure analytical-numerical. Compared against conventional HAM, the technique developed here accelerates the convergence of the solution rapidly. This is important while dealing with problems that cannot be treated beyond second order due to limitations on the computational complexity a typical symbolic algebra software can handle.

Our work comes closest to the optimal homotopy asymptotic method (OHAM). The method has been applied to many problems including non-smooth oscillators \cite{LIAO20102003,NIU20102026,JIA2017865,Fan2013,HERISANU20101607}. As summarized in \cite{HERISANU20101607}, OHAM involves an auxiliary function, similar to our convergence-control function $h(p),$ but more general by assuming it to be the function of both $p$ and time $t$. Authors assume it to be a polynomial in $p$ with unknown coefficients which are to be found out using one of the methods of weighted residual. Convergence crucially depends on the assumed form. Also nonlinear operator in the homotopy involves an arbitrary parameter which when determined using the principle of minimal sensitivity gives essentially the natural frequency of the periodic response. Our approach is more natural as well as simpler. We do not assume any form of the convergence-control function, instead homotopy analysis framework naturally provides the unknowns to be found using the chosen method of weighted residual. We introduce time-stretching function that is assumed to be a function of $p$, which eliminates the need for any arbitrary parameter in the nonlinear operator. We also differ in the important step of the removal of secular terms from their work.

We begin by modelling the oscillations of a spring-mass system with asymmetric restraints using Heaviside function and obtain the approximate periodic solution using Linstedt-Poincaré perturbation method (section \ref{section_2}). We then develop the homotopy analysis and Galerkin projections framework for the same oscillator in section \ref{section_3}. In section \ref{section_4}, we consider periodic solutions of oscillators with non-smoothness of type signum and modulus following the same framework developed in the previous section. Finally in section \ref{sec::5}, we handle the unilaterally constrained simple pendulum modelling it using the Dirac-delta-like function. The decaying response of such oscillator is captured using the homotopy that involves the carefully chosen damped linear oscillator as a linear operator. We compare the natural frequency obtained for each oscillator with the one obtained via methods such as numerical integration (MATLAB ode45), Linstedt-Poincaré perturbation, the method of harmonic balance, averaging combined with non-smooth temporal transform and HAM. The comparison illustrates the usefulness of our framework. Solutions of oscillators discussed in sections \ref{section_3} and \ref{sec::5} are scalable, resulting in the amplitude-independent natural frequency while the same of oscillators in section \ref{section_4} are amplitude dependent. Within the same framework, we are able to approximate periodic as well as non-periodic solutions of the non-smooth oscillators.

\section{A spring-mass system with asymmetric restraints}
\label{section_2}
\subsection{The system}  
\begin{figure}[h]
\centering{\includegraphics[scale=0.4]{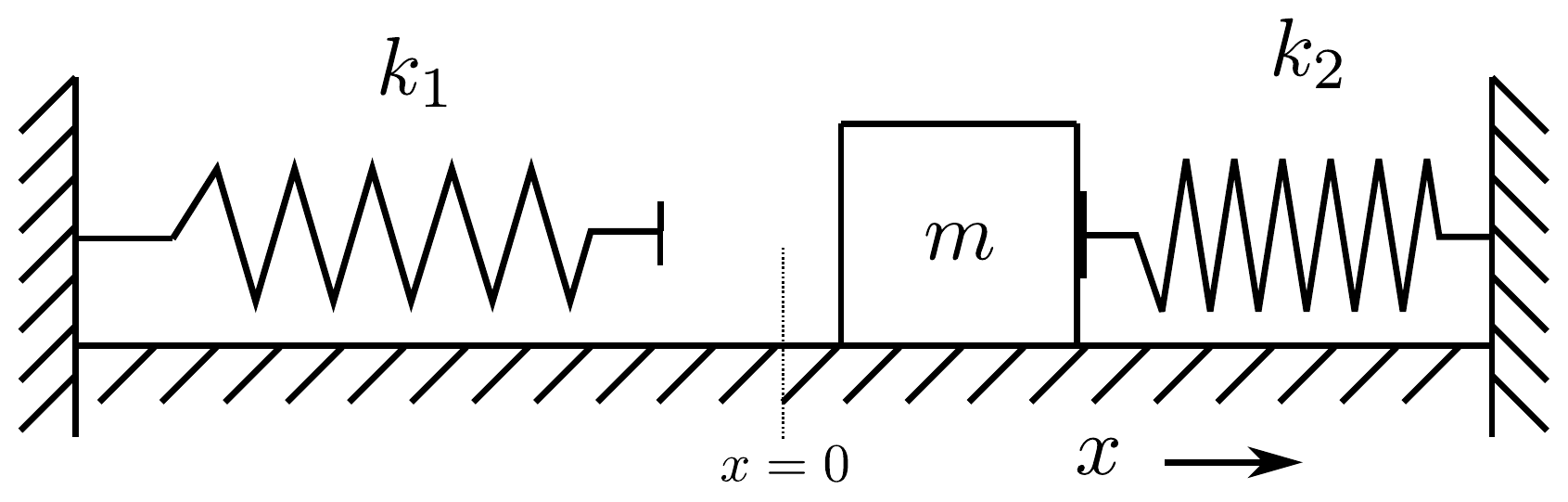}}
\caption{A spring-mass system with non-smoothness of type Heaviside}
\label{springs}
\end{figure} 
Consider mass $m$ kept on a frictionless surface restrained by two springs of stiffness $k_1$ and $k_2$ ($k_1 \neq k_2$) as shown in Fig.\ref{springs}.
The mass vibrates about mean position $x=0$ for some non-zero initial condition. While $x >0,$ the mass loses contact with the spring of stiffness $k_1$ and with the spring of stiffness $k_2$ while $x<0.$ The system stiffness $k$ is given by
\begin{equation*}
   k= \left\{
\begin{array}{ll}
      k_1 \quad x\leq 0 \\
      k_2 \quad x> 0. 
\end{array} 
\right. 
\end{equation*}
Applying the force balance and modelling the system stiffness discontinuity using Heaviside function, we arrive at equation governing the motion
\begin{equation*}
m\ddot{x} + \Big((k_2-k_1) H(x) + k_1\Big)x = 0,
\end{equation*}
where $H$ denotes the Heaviside function. Defining $\epsilon=\tfrac{k_2}{k_1}-1$ and scaling time by $\sqrt{\tfrac{k_1}{m}},$ we simplify the above to
\begin{align}\label{ODE_heaviside}
\ddot{x} +  (1+ \epsilon H(x))x = 0,
\end{align}
where the overdot now denotes differentiation w.r.t. to the scaled time. Equation (\ref{ODE_heaviside}) is a non-smooth oscillator: linearisation near $x=0$ is not possible. Since $H(x_1+x_2)\neq H(x_1) + H(x_2),$ principle of superpositon does not hold; however the solutions are scaleable. Therefore, unlike a nonlinear oscillator, the natural frequency is independent of initial conditions which we attempt to approximate using Lindstedt-Poincaré perturbation method here and using homotopy analysis method combined with Galerkin projections in the next section.

\subsection{Natural response of Eqn.(\ref{ODE_heaviside}) using Lindstedt-Poincaré perturbation method}

\begin{figure}[!h]  
	\begin{subfigure}[b]{0.49\textwidth}
		\centering
		\includegraphics[draft=\status ,height=4cm,width=8cm]{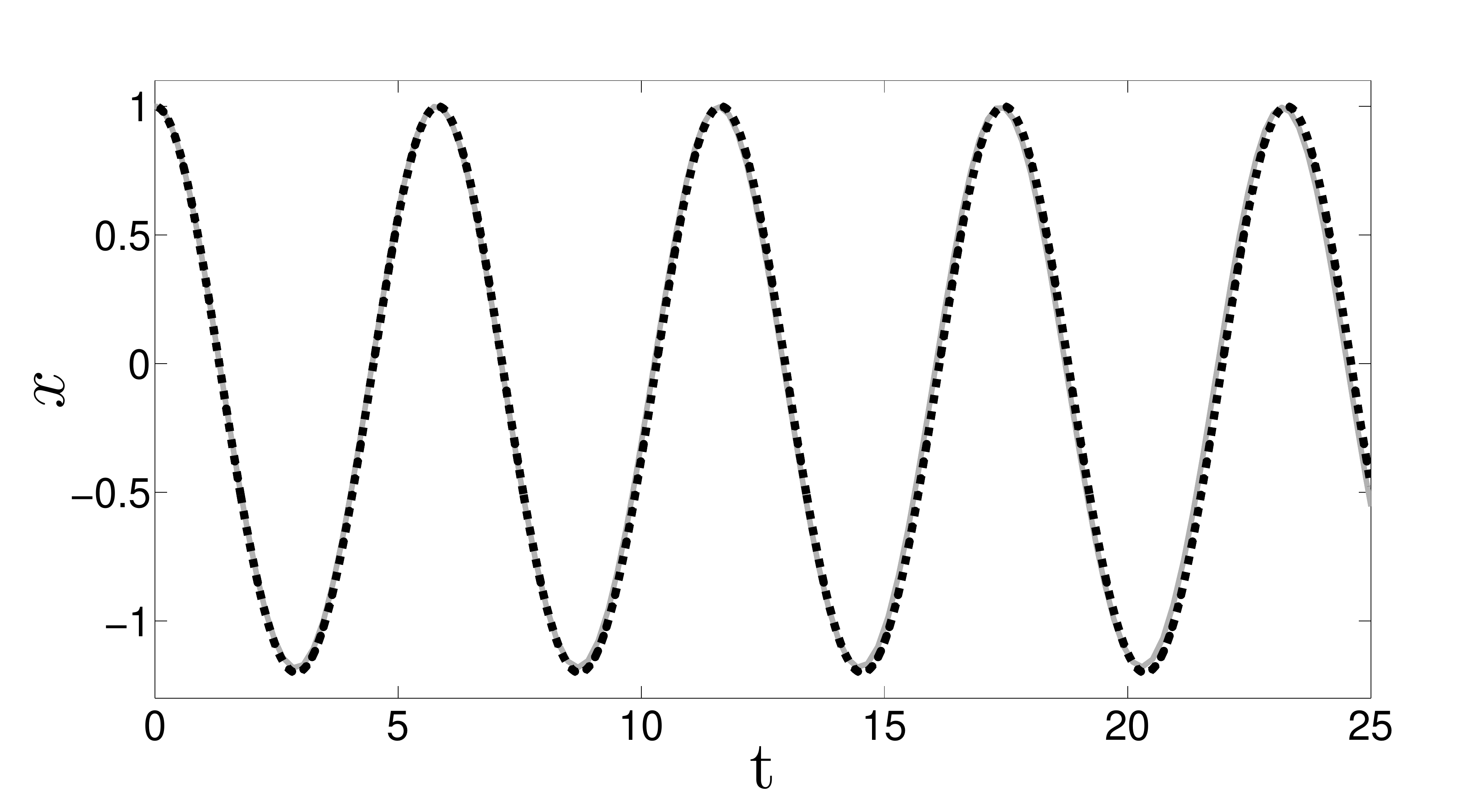}
		\caption{$\epsilon=0.4$}
		\end{subfigure}
	\begin{subfigure}[b]{0.49\textwidth}
  		\centering
		\includegraphics[draft=\status ,height=4cm,width=8cm]{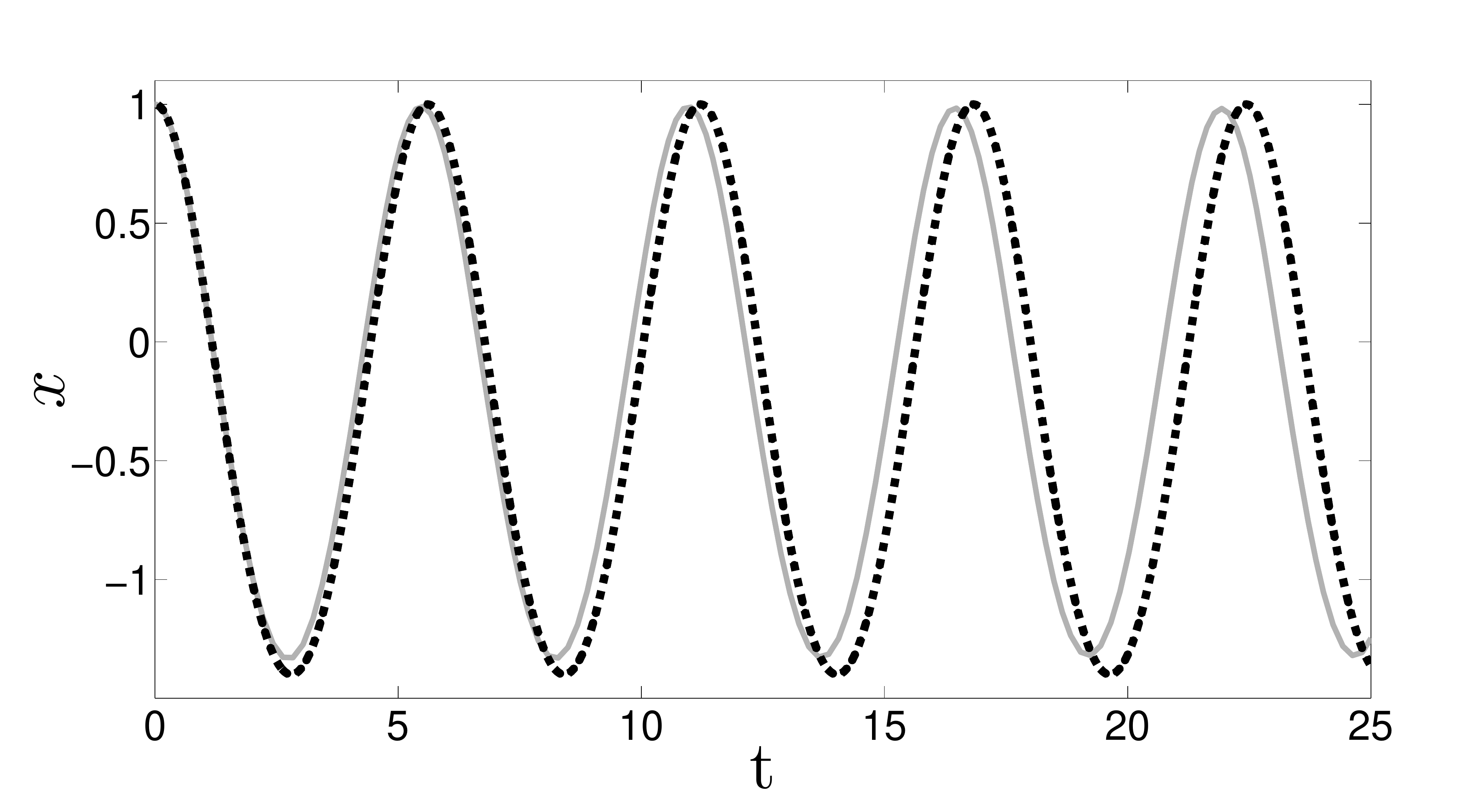}
		\caption{$\epsilon=0.8$}
	\end{subfigure}
	\caption{Comparison of the periodic solution obtained by LP ($x(t)$) and ode$45$; dashed line - LP, continuous line - ode$45.$}
	\label{heaviside_LPP}
\end{figure}

Lindstedt-Poincaré (LP) method is a singular perturbation technique employed for obtaining periodic solutions of weakly nonlinear oscillators. It is unusual to apply LP method to a non-smooth oscillator; however a few attempts have been made \cite{Vestroni2008}. Here, we assume values of spring stiffness $k_1$ and $k_2$ such that the perturbation parameter assumes $0<\epsilon \ll 1.$ We attempt to approximate a periodic solution to Eqn.(\ref{ODE_heaviside}) with initial condition $(1,0)$ as the solutions to Eqn.(\ref{ODE_heaviside}) are scaleable. We begin by scaling time $t$ to $\tau=\Omega t$ with
\begin{equation}\label{taylor_omega_LPP}
\Omega = 1 + \epsilon \omega_1 + \epsilon^2 \omega_2 + \mathcal{O}(\epsilon^3),
\end{equation}
where $\omega_i'$s are unknowns. Equation (\ref{ODE_heaviside}) in the scaled time is
\begin{equation}\label{ODE_heaviside_scaled_time}
\Omega^2 {x}_{\tau\tau} +  (\epsilon H(x) + 1)x = 0,
\end{equation}
where the subscript $\tau$ denotes differentiation w.r.t. $\tau.$ Solution to Eqn.(\ref{ODE_heaviside}) in the scaled time $\tau$ is assumed to be a series strictly asymptotic in $\epsilon$ and is expressed as
\begin{equation}\label{taylor_x_LPP}
x(\tau) = x_0(\tau) + \epsilon x_1(\tau) + \epsilon^2 x_2(\tau) + \mathcal{O}(\epsilon^3).
\end{equation}
Each $x_i(\tau),\;i=0,1,\cdots$ is periodic and $\mathcal{O}(1)$ for all $\tau.$ Substituting Eqn.(\ref{taylor_omega_LPP}) and Eqn.(\ref{taylor_x_LPP}) in Eqn.(\ref{ODE_heaviside_scaled_time}), we get
\begin{multline}\label{ODE_LPP}
x_{0,\tau\tau} + x_0 + \epsilon \Big( x_0 H ( x_0 +\epsilon\,x_1 +{\epsilon}^{2}x_2 ) + x_1 +2 \omega_{{1}}x_{0,\tau\tau} +x_{1,\tau\tau} \Big) +\\ \epsilon^2 \Big( x_1 H ( x_0 +\epsilon x_1 +{\epsilon}^{2}x_2 ) + x_2 + {\omega_{{1}}}^{2}x_{0,\tau\tau} + 2 \omega_{{2}} x_{0,\tau\tau} +2\omega_{{1}}x_{1,\tau\tau}+x_{2,\tau\tau}\Big) + \cdots=0.
\end{multline}
To obtain the coefficients at various orders of $\epsilon$ explicitly from the above and with the assumption of $0<\epsilon\ll 1,$ we  approximate
\begin{equation}\label{heaviside_approx}
H(x_0 +\epsilon\,x_1 +{\epsilon}^{2}x_2+\cdots) \approx H(x_0) + \epsilon x_1\delta(x_0) + \epsilon^2\left(\delta(x_0)x_2 + \frac{1}{2}\delta'(x_0)x_1^2\right). 
\end{equation}
The above approximation may become exact for certain oscillators, for example, Eqn.(\ref{ODE_signum}). Simplifying Eqn.(\ref{ODE_LPP}) using Eqn.(\ref{heaviside_approx}) and then collecting the coefficients of different powers of $\epsilon,$ we get at various orders
\begin{subequations}
\begin{gather}
\mathcal{O}(\epsilon^0): \ x_{0,\tau\tau} + x_{{0}} = 0,\label{LPP_epsilon_0}\\
\mathcal{O}(\epsilon^1): \ x_{1,\tau\tau} + x_1 = -\ x_0 H ( x_0) - 2 \omega_{{1}}x_{0,\tau\tau},\label{LPP_epsilon_1} \\
\mathcal{O}(\epsilon^2): \ x_{2,\tau\tau} + x_2 = -x_1 H ( x_0 ) - {\omega_{{1}}}^{2}x_{0,\tau\tau}  - x_0x_1\delta(x_0)-2 \omega_{{2}} x_{0,\tau\tau} - 2\omega_{{1}}x_{1,\tau\tau}\label{LPP_epsilon_2}.
\end{gather}
\end{subequations}
Solving Eqn.(\ref{LPP_epsilon_0}) for $x_0(\tau)$ with initial condition $(1,0)$, we get
\[x_0(\tau) = \cos\tau. \]
Substituting the above expression for $x_0(\tau)$ in Eqn.(\ref{LPP_epsilon_1}), we obtain
\begin{equation}\label{first_order_LPP}
x_{1,\tau\tau} + x_1 = -\cos\tau H \Big(\cos\tau\Big) + 2 \omega_{{1}}\cos\tau.
\end{equation}
Removal of secular terms, the central aspect of any singular perturbation method is not that straightforward for Eqn.(\ref{first_order_LPP}). Since $x(\tau)$ and hence $x_1(\tau)$ are periodic with period $2\pi,$ it suffices to construct both over one time period,  $\tau\in[0,2\pi].$ By restricting $H(\cos\tau)$ over $\tau\in[0,2\pi],$ the r.h.s. of Eqn.(\ref{first_order_LPP}) maybe modified as
\begin{equation}\label{first_order_LPP_new}
x_{1,\tau\tau} + x_1 = -\cos\tau \bigg(1-H \left(\tau-\frac{\pi}{2}\right) + H \left(\tau-\frac{3\pi}{2}\right)\bigg) + 2 \omega_{{1}}\cos\tau =:F_1(\tau).
\end{equation}
Equation (\ref{first_order_LPP_new}) is a restricted version of Eqn.(\ref{first_order_LPP}) over the interval $\tau\in[0,2\pi].$ It is forced externally at its natural frequency. We eliminate the secular term from the right hand side of Eqn.(\ref{first_order_LPP_new}) using
\begin{equation}\label{secular_term}
\int_0^{2\pi} F_1(\tau) \cos\tau d\tau = 0 .
\end{equation}
Solving Eqn.(\ref{secular_term}) for $\omega_1$ gives
\[ \omega_1 = \frac{1}{4}.\]
Solving Eqn.(\ref{first_order_LPP_new}) for $x_1(\tau)$ with $\omega_1=\tfrac{1}{4}$ and initial condition $(0,0),$  we get
\begin{equation*}
x_1(\tau) = \left( \frac{1}{2}\cos\tau + \left( -\frac{\pi}{4}+\frac{\tau}{2} \right) \sin\tau  \right) {\it H} \left( \tau-\frac{\pi}{2} \right)
+\left( -\frac{1}{2}\cos\tau + \left( \frac{3\pi}{4} -\frac{\tau}{2} \right) 
\sin\tau  \right) {\it H} \left( \tau-\frac{3\pi}{2}\right) -\frac{\tau}{4}\sin\tau.
\end{equation*}
Though the expression for $x_1(\tau)$ is linearly growing w.r.t. $\tau,$ it is valid only over the finite interval $[0,2\pi].$ Repeating the above steps in case of Eqn.(\ref{LPP_epsilon_2}), we get
\[\omega_2 = -\frac{1}{8}\quad \mbox{and}\]
\begin{dmath*}
x_2(\tau)=
 \left(  \left( -\frac{{\pi }^{2}}{32}+\frac{\pi \tau}{16}-\frac{1}{8} \right) \cos\tau + \left( \frac{\pi}{16}-\frac{\tau}{8} \right) \sin\tau 
 \right) H \left( \tau-\frac{\pi}{2} \right) +\left(  \left( -{\frac {3\,{\pi }^{2}}{32}}+\frac{\pi \tau}{16}+\frac{1}{8}\right) \cos\tau + \left( -\frac{3\pi}{16} + \frac{\tau}{8} \right) \sin\tau  \right) H \left( \tau-\frac{3\pi}{2} \right) 
 -\frac{{\tau}^{2}}{32}\cos\tau +\frac{\tau}{16}\sin\tau.
\end{dmath*} 
Substituting for $\omega_1$ and $\omega_2$ in Eqn.(\ref{taylor_omega_LPP}), we get
\begin{equation}\label{omega_LPP}
\Omega = 1 + \frac{\epsilon}{4} - \frac{\epsilon^2}{8}=:\omega_{LP}.
\end{equation}
Finally, substituting the expressions for $x_0(\tau),\ x_1(\tau)$ and $x_2(\tau)$ in Eqn.(\ref{taylor_x_LPP}) and re-scaling time back to $t,$ we obtain one period of the periodic solution. Figure (\ref{heaviside_LPP}) shows the solution obtained via LP as well as via MATLAB numerical integrator ode$45$ for two values of $\epsilon.$ It can be seen clearly that as $\epsilon$ increases, the match between the two solutions deteriorates. One of the major reasons for this mismatch, apart from the loss of asymptotic nature of the series (Eqn.\eqref{taylor_x_LPP}) is the approximation  Eqn.(\ref{heaviside_approx}).

\section{Natural response of Eqn.\ref{ODE_heaviside} using homotopy analysis and Galerkin projections}
\label{section_3}
\label{HG:heaviside}
Smallness of the perturbation parameter $\epsilon$ that plays a crucial role in the Linstedt-Poincaré method is not a restriction as far as the application of homotopy analysis method (HAM) to Eqn.(\ref{ODE_heaviside}) is considered. We again attempt a periodic solution $x(t)$ to Eqn.(\ref{ODE_heaviside}) with initial condition $(1,0)$ and begin by introducing an embedding parameter $p \in [0,1],$ deformation $\tilde{x}(t;p)$ and construct homotopy $\mathcal{H}$ as
\begin{equation}\label{PDE_heaviside_unscaled}
	\mathcal{H} \equiv (1-p)\mathcal{L}\Big(\tilde{x}(t;p) \Big) - h(p) \mathcal{N}\Big(\tilde{x}(t;p) \Big) = \ 0.
\end{equation}
Here, we choose $h(0)=0$ and linear operator $\mathcal{L}\equiv \ddot{x}+x.$ Nonlinear operator $\mathcal{N}$ is set to  Eqn.(\ref{ODE_heaviside}). As $p$ varies continuously from $0$ to $1,$ the deformation $\tilde{x}(t;p)$ varies continuously from the solution of $\mathcal{L}\Big(\tilde{x}(t;0)\Big)=0$ to the solution of  $\mathcal{N}\Big(\tilde{x}(t;1)\Big)=0.$ $\tilde{x}(t;1)$ is the desired periodic solution $x(t).$ Usual application of homotopy analysis considers the convergence-control parameter $h$ as a constant. We introduce it in Eqn.(\ref{PDE_heaviside_unscaled}) as a function of the embedding parameter $p$ and hence call it a convergence-control function. Making it a function of $p$ introduces extra unknowns (its derivatives at $p=0$) which are to be found using extra equations via Galerkin projections. We subject it to conditions $h(0)=0$ and $h(1) \neq 0.$

We scale time $t$ to $ \tau=\omega(p)t $ and introduce a time-stretching function $\lambda(p)=\tfrac{1}{\omega(p)}.$ We also choose $\lambda(0)=1.$ The frequency of the solution to Eqn(\ref{ODE_heaviside}) in this framework is $\omega(1)=\frac{1}{\lambda(1)}.$ Homotopy with respect to the scaled time is
\begin{equation}\label{PDE_heaviside_scaled}
\mathcal{H} \equiv (1-p) \Big(\tilde{x}_{\tau\tau} + \lambda(p)^2 \tilde{x}\Big) -h(p)\Big(\tilde{x}_{\tau\tau} + \lambda(p)^2(\epsilon H(\tilde{x}) + 1)\tilde{x}  \Big)=0.
\end{equation}
Taylor-expanding $\tilde{x}(\tau;p)$ about $p=0,$ we have
\begin{align}
	\tilde{x}(\tau;p) = x^{\mbox{\tiny \mbox{\tiny [0]}}}(\tau) + \sum_{n=1}^{\infty}\frac{1}{n!} x^{\mbox{\tiny [{\em n}]}}(\tau)p^{n},\qquad \qquad \qquad \quad \ \label{solution_taylor}\\
	\mbox{where}  \quad  x^{\mbox{\tiny \mbox{\tiny [0]}}}(\tau) = \tilde{x}(\tau;0)  \quad \mbox{and} \quad  x^{\mbox{\tiny [{\em n}]}}(\tau) = \left. \pd{n}{\tilde{x}}{p} \right|_{p=0}. \nonumber
\end{align}
Here as well, we expect all $x^{\mbox{\tiny [{\em n}]}}$ to be periodic. Differentiating Eqn.(\ref{solution_taylor}) w.r.t. $\tau$, we get
\begin{equation}\label{solution_taylor_xdot}
   \tilde{x}_{\tau} = {x}^{\mbox{\tiny \mbox{\tiny [0]}}}_{\tau}(\tau) + \sum_{n=1}^{\infty}\frac{1}{n!} {x}^{\mbox{\tiny [{\em n}]}}_{\tau}(\tau)p^{n}.  
\end{equation}
We also Taylor-expand the time-stretching and convergence-control functions
\begin{align}
	\lambda(p) = 1 + \sum_{n=1}^{\infty}\frac{1}{n!} \lambda^{\mbox{\tiny [{\em n}]}}p^{n} \quad &\mbox{where} \quad \lambda^{\mbox{\tiny [{\em n}]}}=  \left.\frac{d^{n}\lambda(p)}{dp^{n}} \right|_{p=0}\label{solution_taylor_lambda},\\
	h(p) =  \sum_{n=1}^{\infty}\frac{1}{n!} h^{\mbox{\tiny [{\em n}]}}p^{n} \quad &\mbox{where} \quad h^{\mbox{\tiny [{\em n}]}} =  \left.\frac{d^{n}h(p)}{dp^{n}} \right|_{p=0}.\label{solution_taylor_h}
\end{align}

We begin by obtaining $x^{\mbox{\tiny [0]}}(\tau).$ Substituting $p=0$ in Eqn.(\ref{PDE_heaviside_scaled}), we get the zeroth-order deformation equation
\begin{equation}\label{zeroth_order_eqn_impulse}
{x}^{\mbox{\tiny[0]}}_{\tau\tau}+ x^{\mbox{\tiny [0]}} = 0.
\end{equation}
By choosing the initial conditions for the above oscillator as $(1,0),$ we get
\begin{equation}\label{zeroth_order_solution}
 x^{\mbox{\tiny [0]}}(\tau) = \cos\tau.
\end{equation}
Our choice of the initial condition for zeroth-order deformation equation along with Eqns.(\ref{solution_taylor}) and (\ref{solution_taylor_xdot}) fixes the initial conditions for the higher order deformation equations as
\begin{equation}\label{tmp}
	x^{\mbox{\tiny[{\em n}]}}(0)=0, \quad {x}^{\mbox{\tiny[{\em n}]}}_{\tau}(0)=0, \quad n\geq1. 
\end{equation}
To get the first-order deformation equation, we differentiate Eqn.(\ref{PDE_heaviside_scaled}) w.r.t. $p$ once. Then substituting $p=0,$ and using Eqns.\eqref{solution_taylor},\eqref{solution_taylor_lambda},\eqref{solution_taylor_h} and \eqref{zeroth_order_solution}, we get
\begin{dmath*}
 {x}^{\mbox{\tiny [1]}}_{\tau\tau} + {x}^{\mbox{\tiny [1]}} = 
\bigg(H(\cos\tau) \epsilon h^{\mbox{\tiny[1]}} -2 \lambda^{\mbox{\tiny \mbox{\tiny [1]}}} \bigg) \cos\tau.
\end{dmath*}
As we did in the previous section, we simplify the above equation by restricting it to the interval of interest, i.e., $\tau=[0,2\pi].$ The first-order deformation equation is then modified to
\begin{equation}\label{first_order_new}
 {x}^{\mbox{\tiny [1]}}_{\tau\tau} + {x}^{\mbox{\tiny [1]}} = 
\bigg(\bigg(1-H \left(\tau-\frac{\pi}{2}\right) + H \left(\tau-\frac{3\pi}{2}\right)\bigg)\epsilon h^{\mbox{\tiny[1]}} - 2 \lambda^{\mbox{\tiny \mbox{\tiny [1]}}} \bigg) \cos\tau =: F_1(\tau).
\end{equation}
Removing the secular term using Eqn.(\ref{secular_term}), we fetch
\[	\lambda^{\mbox{\tiny \mbox{\tiny [1]}}} = \frac{\epsilon h^{\mbox{\tiny[1]}} }{4}.
\]
We solve Eqn.(\ref{first_order_new}) using initial condition from Eqn.(\ref{tmp}) to get
\begin{equation}\label{first_order_solution}
{x}^{\mbox{\tiny [1]}}(\tau) = \epsilon h^{\mbox{\tiny[1]}} \Bigg(\frac{\tau}{4}\sin\tau + \bigg( \frac{1}{2}\cos \tau + \left( -\frac{3\pi  }{4} +\frac{\tau}{2} \right) \sin\tau\bigg)H \left( \tau-\frac{3\pi}{2}\right) + \left( 
-\frac{1}{2}\cos\tau + \left( \frac{\pi}{4} -\frac{\tau  }{2} \right) \sin\tau \right)H\left( \tau-\frac{\pi}{2}\right)\Bigg).
\end{equation}
Following the same procedure, we proceed upto second order and eliminate the secular term to obtain
\begin{equation*}
	\lambda^{\mbox{\tiny [2]}} = \frac{\epsilon}{2}\left({h^{\mbox{\tiny[1]}}}+\frac{ h^{\mbox{\tiny[2]}}}{2}+{h^{\mbox{\tiny[1]}}}^{2}\right)+\frac{3{\epsilon}^{2}{h^{\mbox{\tiny[1]}}}^{2}}{8}.
\end{equation*}
Substituting $p=1$ in Eqn.(\ref{solution_taylor_lambda}) and expressions for $\lambda^{\mbox{\tiny [1]}},\ \lambda^{\mbox{\tiny [2]}}$ from the above, we get
\begin{equation*}
	\lambda(1) = 1 + \frac{\epsilon}{2}\left(\frac{3{h^{\mbox{\tiny[1]}}}}{2}+\frac{ h^{\mbox{\tiny[2]}}}{2}+{h^{\mbox{\tiny[1]}}}^{2}\right)+\frac{3{\epsilon}^{2}{h^{\mbox{\tiny[1]}}}^{2}}{8}+\cdots.
\end{equation*}
The frequency of the desired solution $\omega(1)$ is then given by
\begin{equation}\label{omega_HAM}
	\omega(1)=\frac{1}{\lambda(1)} \approx 1 - \frac{\epsilon}{2}\left(\frac{3{h^{\mbox{\tiny[1]}}}}{2}+\frac{ h^{\mbox{\tiny[2]}}}{2}+{h^{\mbox{\tiny[1]}}}^{2}\right) + {\epsilon}^{2}\left(\frac{1}{4}\left(\frac{3{h^{\mbox{\tiny[1]}}}}{2}+\frac{ h^{\mbox{\tiny[2]}}}{2}+{h^{\mbox{\tiny[1]}}}^{2}\right)^2-\frac{3{h^{\mbox{\tiny[1]}}}^{2}}{8}\right) =: \omega_{HG}.
\end{equation}
Ignoring higher order terms, substituting $p=1$ in Eqn.(\ref{solution_taylor}), using Eqns.(\ref{zeroth_order_solution}) and (\ref{first_order_solution}) and scaling time back to $t,$ we get
\begin{equation}\label{taylor_approx4} 
	\tilde{x}(t;1) \approx x^{\mbox{\tiny[0]}}(t) + x^{\mbox{\tiny[1]}}(t) =: x_1(t).
\end{equation}
Subscript $1$ here denotes that we approximate solution by solving upto first-order deformation equation. The solution $x_1(t)$ as obtained using Eqn.(\ref{taylor_approx4}) is restricted to $t \in [0,\tfrac{2\pi}{\omega(1)}].$ Since we approximate one period of the periodic solution, we have the desired solution for all time $t.$ The resulting $x_1(t)$ is not only a function of $\omega(1)$ but also of free parameters $h^{\mbox{\tiny[1]}}$ and $h^{\mbox{\tiny[2]}}.$ Therefore, we need more equations in addition to Eqn.(\ref{omega_HAM}) involving these unknowns. We use Galerkin method of the weighted residual to obtain these equations. Using the expression for $x_1(t)$ from Eqn.(\ref{taylor_approx4}), we get the residual
\begin{equation*}
 \mathcal{R}_1(t) \equiv  \ddot{x}_1(t) +  x_1(t)(\epsilon H(x_1(t)) + 1).
 \end{equation*}
As per the Galerkin scheme, taking two weighting functions $w_1(t)=\cos\Big(\omega(1){t} \Big)$ and $w_2(t)=H\Big( \omega(1)t - \tfrac{\pi}{2} \Big),$ we have
\begin{equation*}
	\int_{0}^{\tfrac{2\pi}{\omega(1)}} w_i(t)\mathcal{R}_1(t) dt = 0, \quad i=1,2.
\end{equation*}
To carry out the above integration, we simplify the term, $H(x_1(t)).$ One straightforward simplification is to replace the argument of Heaviside by the solution to the zeroth-order deformation equation in time $t,$ i.e.,  $x^{\mbox{\tiny[0]}}(t).$ We propose better simplification by finding out zeros of $x_1(t)$ over the interval $[0,\tfrac{2\pi}{\omega(1)}]$ using one iteration of the Newton-Raphson scheme. Denoting the two zeros by $t^{\mbox{\tiny[1]}}$ and $t^{\mbox{\tiny[2]}}$ given by
\begin{equation}
t^{\mbox{\tiny[1]}} = \frac{\pi}{2\omega(1)} - \left.\frac{x_1(t)}{\dot{x}_1(t)}\right|_{t=\frac{\pi}{2\omega(1)}} \quad \mbox{and} \quad t^{\mbox{\tiny[2]}} = \frac{3\pi}{2\omega(1)} - \left.\frac{x_1(t)}{\dot{x}_1(t)}\right|_{t=\frac{3\pi}{2\omega(1)}}, \nonumber 
 \end{equation}
 we approximate $H(x_1(t))$ by,
 \begin{gather}
H(x_1(t)) \approx 1-H(t-t^{\mbox{\tiny[1]}})+ H(t-t^{\mbox{\tiny[2]}}). \label{heaviside_approx_GP}
\end{gather}
We assume that the expressions for $t^{\mbox{\tiny[1]}}$ and $t^{\mbox{\tiny[2]}},$ which are numerically unknown, lie close to $\tfrac{\pi}{2\omega(1)}$ and $\tfrac{3\pi}{2\omega(1)}$ respectively. This results in $H(x_1(t))=0$  between $t^{\mbox{\tiny[1]}}$ and $t^{\mbox{\tiny[2]}}$ and $H(x_1(t))=1$ everywhere else over $[0,\tfrac{2\pi}{\omega(1)}].$ The weighted residual with the approximate expression for $H(x_1(t))$ can now be integrated. Note that all terms in $\mathcal{R}_1(t)$ are of type $\sin(\omega(1)t)$ or $\cos(\omega(1)t)$ or similar. Therefore, after putting the limits of integration, unknown $\omega(1)$ is eliminated and results in $-\cos(2\pi)$ or $\sin(2\pi)$ or similar. The integral so obtained has Heaviside terms with arguments involving $\omega(1),\ h^{\mbox{\tiny[1]}}$ and $h^{\mbox{\tiny[2]}}.$ Fixing $\epsilon,$ we solve the system of above two equations and Eqn.(\ref{omega_HAM}) simultaneously for the unknowns using fsolve command of symbolic algebra package Maple. In order to avoid unwanted roots and reduce the computation burden, we fix the interval for $\omega(1)$ in the neighbourhood of $1$ or $-1.$ 
\vspace{0.1in}

Since the solutions to Eqn.(\ref{ODE_heaviside}) are scalable, the frequency of its solution is independent of initial conditions. This maybe verified by applying the homotopy analysis and Galerkin projections (HG) framework starting with an arbitrarily chosen initial condition instead of $(1,0).$
\begin{figure}[!h]
  \begin{subfigure}[b]{0.49\textwidth}
  \centering
    \includegraphics[draft=\status ,height=4cm,width=8cm]{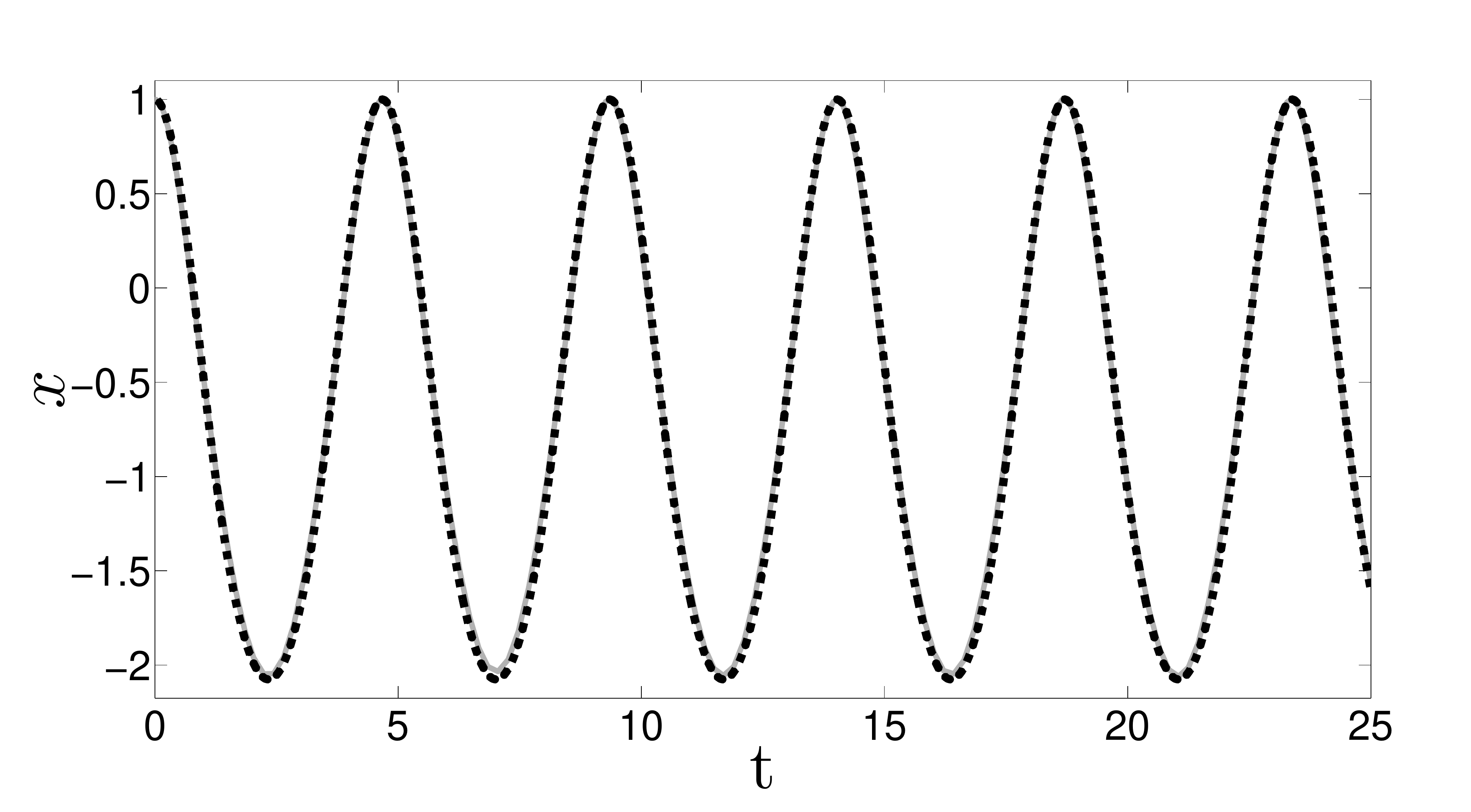}
    \caption{$\epsilon=3.2$}
  \end{subfigure}
  \begin{subfigure}[b]{0.49\textwidth}
\centering
    \includegraphics[draft=\status ,height=4cm,width=8cm]{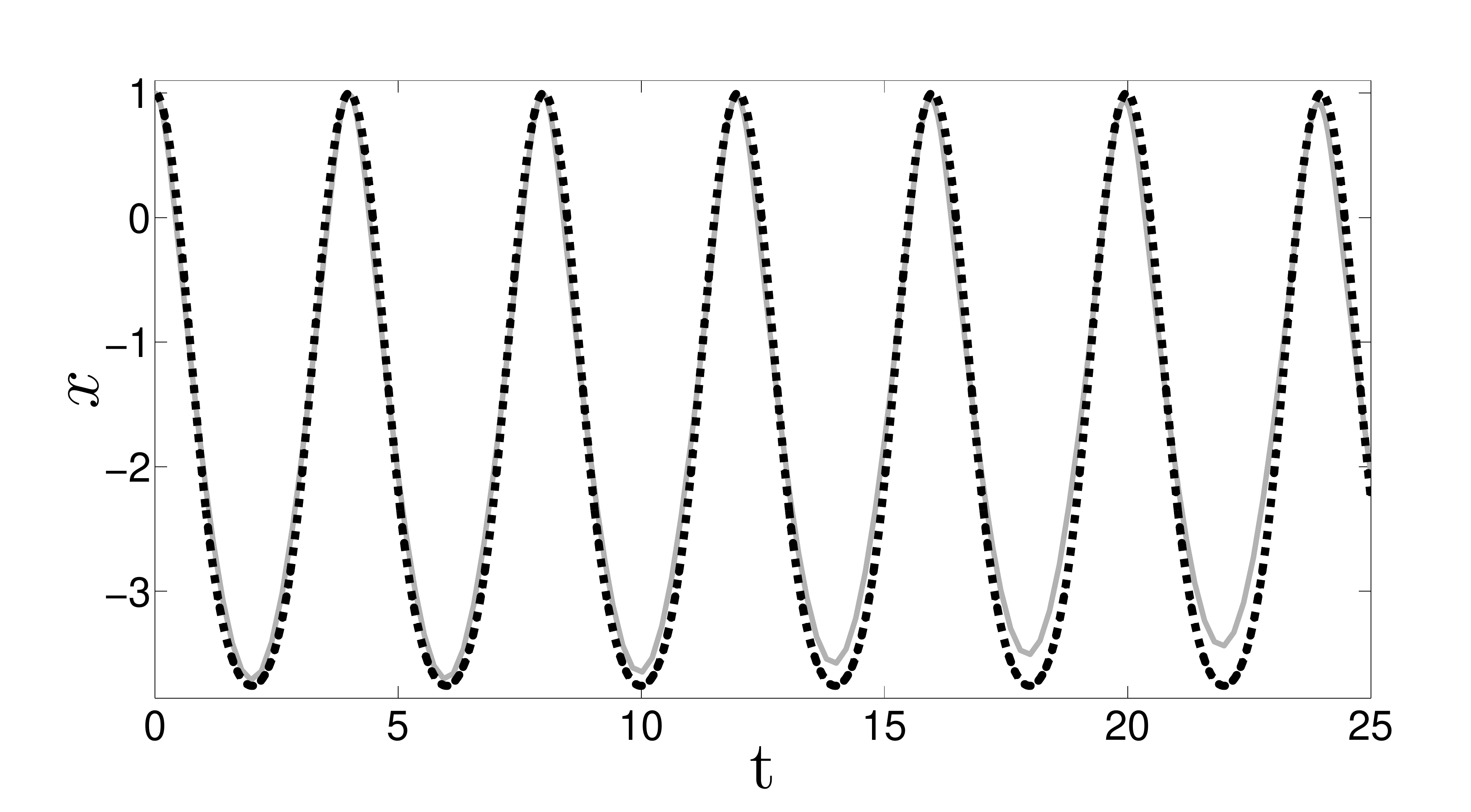}
    \caption{$\epsilon=12.8$}
  \end{subfigure}
  \caption{Comparison of periodic solution obtained by first order homotopy and Galerkin  with ode$45;$ continuous line: ode$45$, dashed line: HG.}
  \label{heaviside_HG}
\end{figure}

We find $x_1(t)$ for two values of $\epsilon$ and compare it with the solution obtained using ode$45$ (built-in integrator of MATLAB) in Fig(\ref{heaviside_HG}). We see that the solution by homotopy analysis and Galerkin projections (HG) mimics the numerical solution for $\epsilon$ as large as $3.2.$ As $\epsilon$ increases, we observe the mismatch between the two solutions in terms of both frequency and amplitude, noticable at $\epsilon=12.8.$ This can be attributed to a number of factors: a significant contribution to the error is due to the approximation of $H(x_1(t))$ (Eqn.\ref{heaviside_approx_GP}). This error may be minimized by either finding a better alternative to Eqn.(\ref{heaviside_approx_GP}), for example, by applying the Newton-Raphson iteration more than once. The choice of weighting functions also plays a significant role in minimizing the error. Here, we set the weighting functions as the basis functions of the solution from HG with coefficients that seem to contribute maximum to the solution amplitude. Lastly, the method is carried out only upto the first order; however at higher orders, the computational burden increases significantly.
In the LP framework, since the series solution (Eqn.(\ref{taylor_x_LPP})) is asymptotic in $\epsilon,$ the validity of the method is restricted to smaller values of $\epsilon.$ For HG method, the Taylor expansion (Eqn.(\ref{solution_taylor})) involves the free parameters that are  determined using a heuristic approach of Galerkin projections. As a result Eqn.(\ref{solution_taylor}) is not guaranteed to be asymptotic in $\epsilon.$
\vspace{0.1in}

Although application of Galerkin projections in HG demands the approximation Eqn\eqref{heaviside_approx_GP}, the contribution is not as significant. One of the major differences between the two techniques is the nature of the expression for frequency. LP yields an expression that is purely analytic. HAM gives an analytical expression involving the free parameters which must be determined numerically using Galerkin projections, thereby making it semi-analytic. Also, the HG method consequently becomes computationally more expensive. Both techniques follow the same treatment of differential equation at every order, i.e., the approximation of $H(\cos\tau)$ and the elimination of secular terms.
\vspace{0.1in}
\begin{figure*}[!h]
\centering{\includegraphics[scale=0.8]{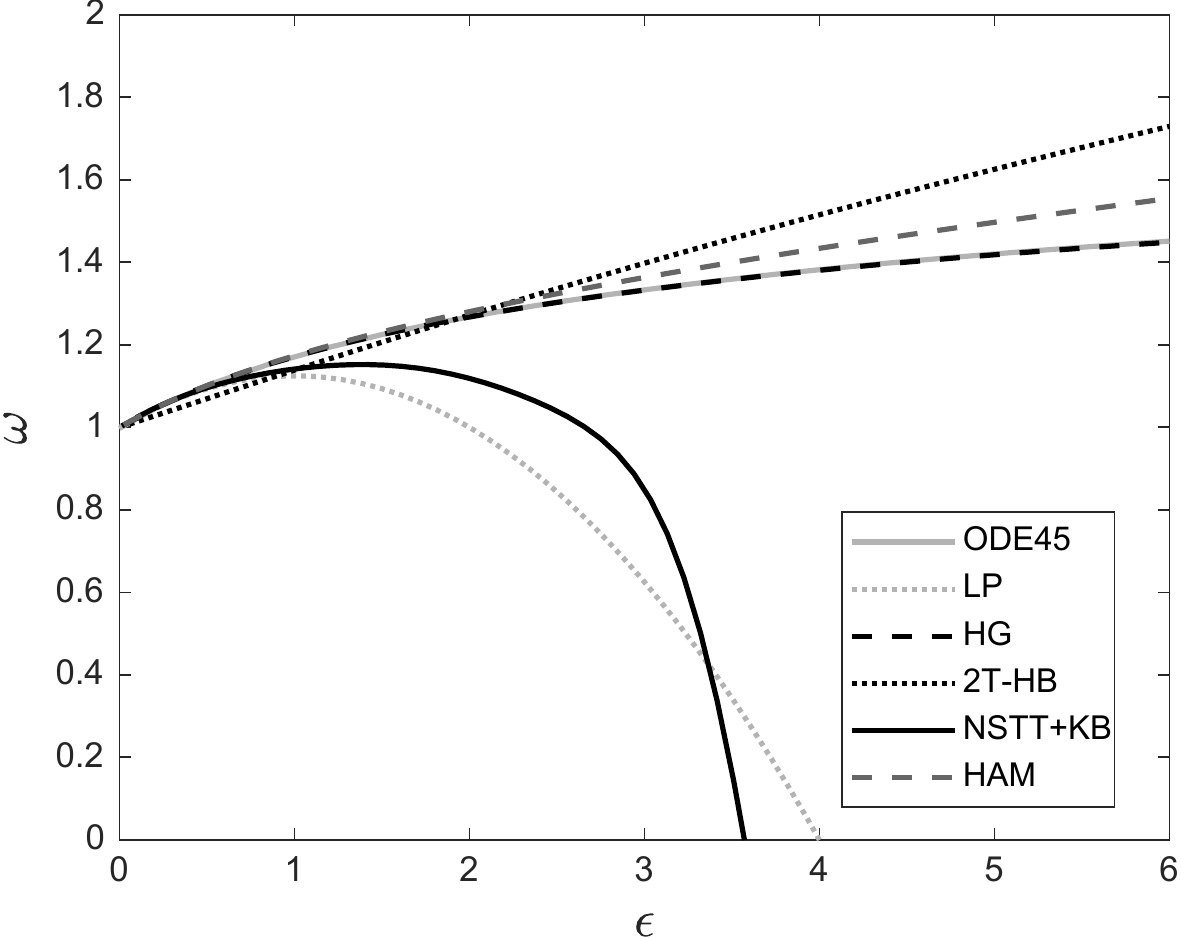}}
\caption{Comparison of the natural frequency $\omega$ obtained via six different methods.}
\label{freq_compare}
\end{figure*} 

We consider the oscillator \eqref{ODE_heaviside} written in the amplitude-phase co-ordinates $[A(t),\varphi(t)]$ as
\begin{equation}\label{ODE_NSTT}
    \begin{aligned}
    \dot{A}&=\frac{\epsilon A}{2}H(A\cos\varphi)\sin 2\varphi,\\
    \dot{\varphi} &= 1-H(A\cos\varphi)\cos^2\varphi.
\end{aligned}
\end{equation}
Here a non-smooth temporal transformation (NSTT) can be applied to $\varphi(t)$ using a triangular sine wave $\tau(2\varphi/\pi)=(2/\pi)\arcsin{\sin (\varphi)}$ and a rectangular cosine wave  $e(2\varphi/\pi)=\mbox{sgn}(\cos(\varphi)),$ transforming Eqn.(\ref{ODE_NSTT}) to
\begin{equation}\label{ODE_NSTT_post}
\begin{aligned}
    \dot{A}&=\frac{\epsilon A}{4}(1+e)\sin(\pi \tau),\\
    \dot{\varphi}&= 1+\frac{\epsilon (1+e)}{2}\cos^2(\pi\tau/2).
\end{aligned}    
\end{equation}
Proposing the above NSTT transformation, the author in \cite{Pilipchuk2015} subjects Eqn.(\ref{ODE_NSTT_post}) to Krylov-Bogolyubov averaging, thus leading to a closed-form analytical  expression for the phase and amplitude computed up to $ \mathcal{O}(\epsilon).$ The resultant solution attains steady state after a few cycles, with a frequency we plot in Fig.(\ref{freq_compare}).

\vspace{0.1in}

 In the context of conventional HAM, authors  prove that the series for natural frequency converges uniquely to a limit independent of the convergence-control parameter $h$ \cite{BPL}. When we carry out the method computationally, the same series will consist of only finitely many terms, and the expression arrived at will depend on $h.$ To select an appropriate value of $h,$ we plot it against the natural frequency for given $\epsilon$ value and the initial condition, and then find the $h$-range where the natural frequency does not change appreciably w.r.t. $h$ \cite{SAPM:SAPM387}. From such a flat region, an arbitrarily chosen value of $h$ gives the approximate natural frequency close to the frequency obtained by computation up to an infinite order.

\begin{figure}[!h]
{\centering{\includegraphics[height=6cm,width=8cm]{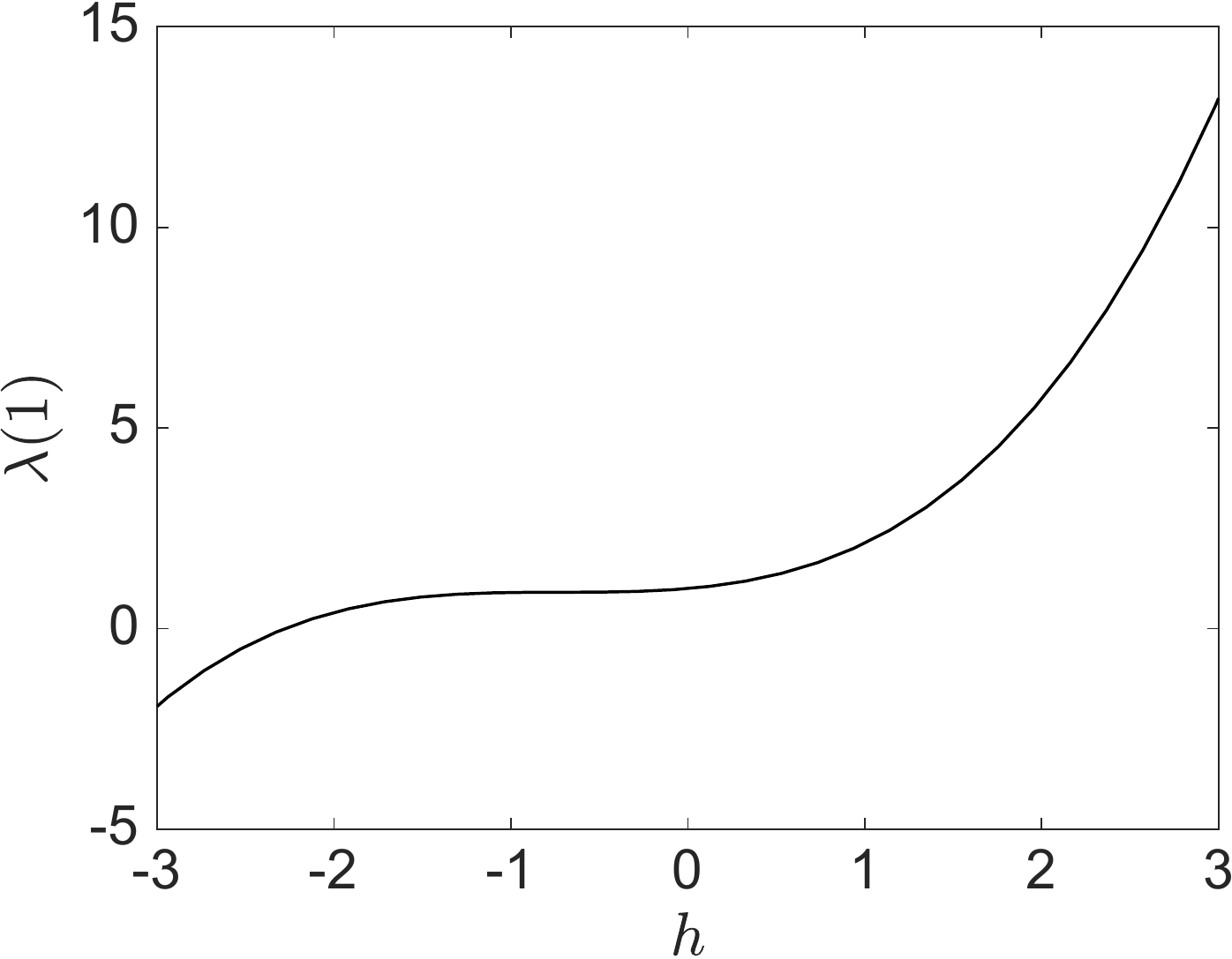}}
\caption{Inverse of frequency $1/\omega(1)$ for $\epsilon=0.5.$}
\label{h_curve}}
\end{figure}

To demonstrate this, we apply HAM to Eqn.(\ref{ODE_heaviside}) by formulating the homotopy that considers $h$ being not a function of $p$,
$$
\mathcal{H} \equiv (1-p)\mathcal{L}\Big(\tilde{x}(t;p) \Big) - h \mathcal{N}\Big(\tilde{x}(t;p) \Big) = \ 0,
$$
with $\mathcal{L}$ and $\mathcal{N}$ as defined in Eqn.(\ref{PDE_heaviside_unscaled}). We take the initial condition as $(1,0),$ scale time $t$ to $\tau=t/\lambda(p)=\omega(p)t,$ compute up to $3$-rd order deformation equation and obtain the expression
\begin{dmath*}
\lambda(1)=1+ \frac{\epsilon}{4} \left( 3h+3{h}^{2}+{h}^{3} \right) +\frac{{\epsilon}^{2}}{8}\left( {\frac {9\,{h}^{2}}{
2}}+3{h}^{3}\right) + {\epsilon}^{3}{\frac {5\,{h}^{3}}{32}}.
\end{dmath*}
Plotting $\lambda(1)$ ($=1/\omega(1)$) against $h$ for $\epsilon=0.5$ as shown in Fig.(\ref{h_curve}), we identify the flat region as $h\in[-1,0].$ In order to select an appropriate $h,$ we use the same argument as mentioned in the previous paragraph for $\lambda(1)$ (instead of frequency) and choose $h=-0.5.$ Carrying the whole procedure out over $\epsilon\in[0,6]$ and choosing appropriate $h$ value every time, we obtain and plot the natural frequency curve via conventional HAM Fig.(\ref{freq_compare}).
\vspace{0.1in}

\vspace{0.1in}

The method of harmonic balance is a special case of Galerkin projections. We also obtain the natural frequency w.r.t. $\epsilon$ using two-term harmonic balance (2T-HB) to Eqn.(\ref{ODE_heaviside_scaled_time}) by taking the ansatz as $A_1\cos \omega t+A_2\cos 3\omega t.$ To determine the unknowns, we compute the residual, generate equations by fixing the initial condition to $(1,0)$ and collect the coefficients of $1,$ $\cos 2t,$ $\cdots$ from the residual. The resultant system of algebraic equations has multiple roots. To determine the frequency for different values of $\epsilon\in[0,5],$ we search for the same in a small neighbourhood of $1$ and plot the results in Fig.(\ref{freq_compare}). Alternatively, we can also avoid multiple roots using the method of generalized harmonic balance method \cite{BELENDEZ20092117}. 
\vspace{0.1in}

Figure (\ref{freq_compare}) compares the frequencies obtained via ode$45,$ LP, NSTT+KB, HAM and $2$T-HB and HG over $\epsilon\in[0,6]$. Treating the frequency obtained via ode$45$ as the most exact, we determine the \% error in the frequencies given by all the above mentioned techniques. We define $\epsilon$-usability range as the range of parameter $\epsilon$ over which the approximate frequency obtained from a given method is in agreement with the same obtained via ode45 within $2\%.$  This range for LP is $\epsilon\in[0,0.8],$ for NSTT+KB  $\epsilon\in[0,1],$ for $3$-rd order HAM $\epsilon\in[0,3],$ for 2T-HB $\epsilon\in[0,2]$ and for HG $\epsilon\in[0,6].$ This establishes the superiority of the hybrid approach of HG over purely asymptotic as well as purely heuristic approaches, though in the limited context of dealing with non-smooth oscillators. Evidently, we observe that NSTT+KB applied up to $\mathcal{O}(\epsilon)$ fares better than LP applied up to $\mathcal{O}(\epsilon^2).$ Also, $2$-nd order HG offers $\epsilon$-usability range that is twice that of $3$-rd order HAM. Using a suitable method of weighted residual to solve for the higher number of unknowns by assuming convergence-control as a function of the embedding parameter, ensures the tighter control over the error leading to more accurate estimate of the natural frequency compared against conventional HAM, thus establishing the accelerated convergence effectively.

\section{Natural response of oscillators involving non-smoothness of type signum and modulus}
\label{section_4}
\subsection{Signum type discontinuity}
\begin{figure}[!h]
	\centering
	  {\includegraphics[draft=\status ,height=4cm,width=8cm]{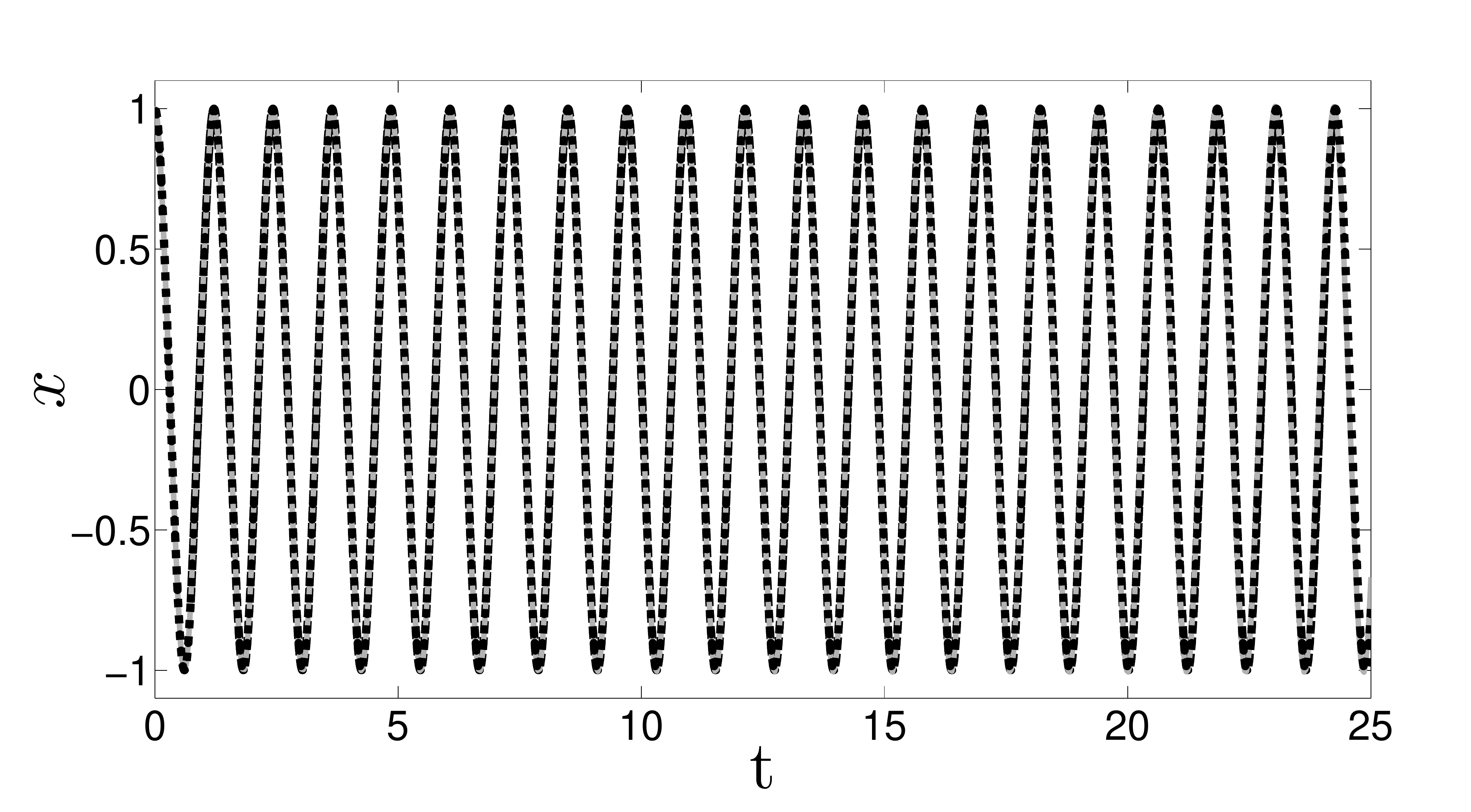}}
	  \caption{Comparison of periodic solution obtained by HG and ode$45$  for $A=1,$ $\epsilon=32$; continuous line: ode$45$, dashed line: HG.}
	  \label{signum_HG}
\end{figure}

\begin{figure*}[h]
\centering{\includegraphics[scale=0.33]{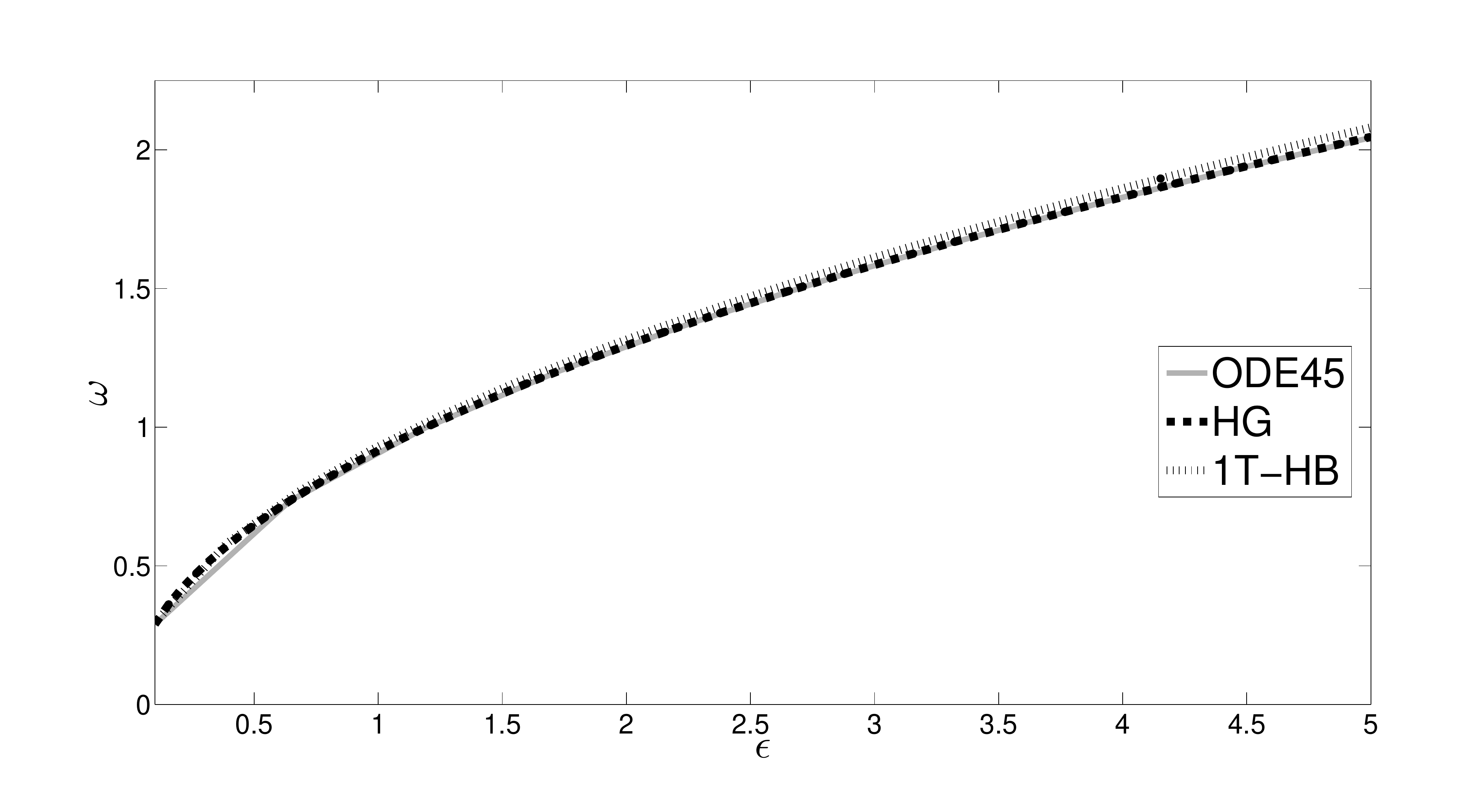}}
\caption{Natural frequency $\omega$ obtained for $A=1$ for oscillator, Eqn.(\ref{ODE_signum_simplified}).}
\label{freq_compare_sig}
\end{figure*}

We consider an oscillator
\begin{equation}\label{ODE_signum}
\ddot{x} + \epsilon x^2 \mbox{sgn}(x) = 0,
\end{equation}
where $\mbox{sgn}$ denotes the signum function. Re-writing the above equation using the Heaviside function, we have 
\begin{equation}\label{ODE_signum_simplified}
\ddot{x} + \epsilon x^2 \Big(2H(x)-1\Big) = 0.
\end{equation}
The above is now treated using HG framework developed in Sec.(\ref{HG:heaviside}). Choosing the linear operator as a simple harmonic oscillator and nonlinear operator as Eqn.(\ref{ODE_signum_simplified}), we construct the homotopy $\mathcal{H}$ and scaling time from $t$ to $\tau=\tfrac{t}{\lambda(p)},$ we get
\begin{equation*}\label{PDE_heaviside_scaled_signum}
\mathcal{H} \equiv (1-p) \Big(\tilde{x}_{\tau\tau} + \lambda(p)^2 \tilde{x}\Big) -h(p)\bigg(\tilde{x}_{\tau\tau} + \epsilon\lambda(p)^2 \tilde{x}^2\Big(2 H(\tilde{x}) - 1\Big) \bigg)=0.
\end{equation*}
The solutions to Eqn.(\ref{ODE_signum}) are not scalable. The oscillator is autonomous though and thereby choosing initial phase to be zero, we attempt a periodic solution with initial condition $(A,0).$ Computing upto the first order, we get approximate expressions for the natural frequency and the solution as 
\begin{multline*}
\omega(1) = 1+h^{\mbox{\tiny[1]}}+\frac{h^{\mbox{\tiny[2]}}}{4}+\frac{{h^{\mbox{\tiny[1]}}}^{2}}{8}-\frac{2\epsilon A}{3\pi}\left( {4{h^{\mbox{\tiny[1]}}}}+{{h^{\mbox{\tiny[2]}}}} +{{{h^{\mbox{\tiny[1]}}}^{2}}}\right) \\-\epsilon^{2}A^2\left(\frac{ \left( 30\,\pi -32 \right) {h^{\mbox{\tiny[1]}}}^{2}}{72\pi}+\frac{4\left( {{4h^{\mbox{\tiny[1]}}}}+{{h^{\mbox{\tiny[2]}}}} +{{{h^{\mbox{\tiny[1]}}}^{2}}}\right)^2}{9\pi^2}\right),
\end{multline*}
\begin{dmath*}
x_1(\tau) = A\cos\tau+\epsilon A^2 h^{\mbox{\tiny[1]}}\left(\left( \frac{1}{3} \cos 2\tau+ \frac{4}{3} \sin\tau - 1 \right) H \left(\tau-\frac{\pi}{2} \right) + \left( -\frac{1}{3} \cos 2\tau+ \frac{4}{3} \sin \tau + 1 \right)H \left(\tau-\frac{3\pi}{2}\right) - \frac{1}{3}\cos\tau -\frac{1}{6}\cos 2\tau - \frac{4\tau}{3\pi}\sin \tau + \frac{1}{2}\right),
\end{dmath*}
\[\mbox{where} \quad \tau=\omega(1){t}. \]
We determine the free parameters $h^{\mbox{\tiny[1]}}$ and $h^{\mbox{\tiny[2]}}$ using Galerkin projections as before using the weighting functions $\cos \left(\omega(1)t\right)$ and $ H \Big(\omega(1)t - \tfrac{\pi}{2} \Big),$ and then solve a system of equations numerically. These unknowns are computed for two values of $\epsilon$ and are given in table (\ref{table:unknowns:sig_abs}). For the same values of $\epsilon,$ Fig(\ref{signum_HG}) shows the comparison between the solution so obtained and the solution via ode$45.$ 
The above obtained expression for $x_1(t)$ has zero-crossings at $t=\frac{\pi}{2\omega(1)},\ {\frac{3\pi}{2\omega(1)}}.$ Hence, $H(x_1(t))$ can be evaluated without any approximation. 
For $A=1,$ the frequency obtained exhibits low \% error for an exceptionally large range of $\epsilon$ (Table (\ref{table:compare:sig})). Due to finite precision of the numerical solver used, we cannot determine frequencies for higher $\epsilon$ values, and hence we are unable to comment on the upper limit of the $\epsilon$-usability range.
\begin{table}[h]
\caption{Frequency comparison for Eqn.(\ref{ODE_signum})}
\centering 
\renewcommand{\arraystretch}{1.25} 
\begin{tabular}{|c| c c | c| } 
\hline
  & \multicolumn{2}{c}{$\omega(1)$} & \multicolumn{1}{|c|}{\% error}\\
\hline 
$\epsilon$ & ode$45$ & HG & HG \\[0.5ex]
\hline  
$3$ &  1.58368 & 1.58521 &     0.09661\\
$12$ & 3.16736 & 3.17042 &     0.09661\\
$50$ &   6.46694 & 6.47159 &   0.07190\\
$200$ & 12.92821 & 12.94318 &  0.11579\\
$800$ &  25.87322 & 25.88637 & 0.05082\\
$8\times10^{17}$ &   $8.17836\times10^8$ & $8.18598\times10^8$ & 0.09317\\
\hline 
\end{tabular}
\label{table:compare:sig}
\end{table}
Equation (\ref{ODE_signum}) is easily amenable to harmonic balance \cite{Mickens}. Figure \eqref{freq_compare_sig} shows the natural frequency comparison using ode$45,$ HG and 1T-HB. Note that LP is not applicable to Eqn.(\ref{ODE_signum}) since  the unperturbed version is not an oscillator. 

\subsection{Modulus type discontinuity}
\begin{figure}[!h]
\centering
  {
    \includegraphics[draft=\status ,height=4cm,width=8cm]{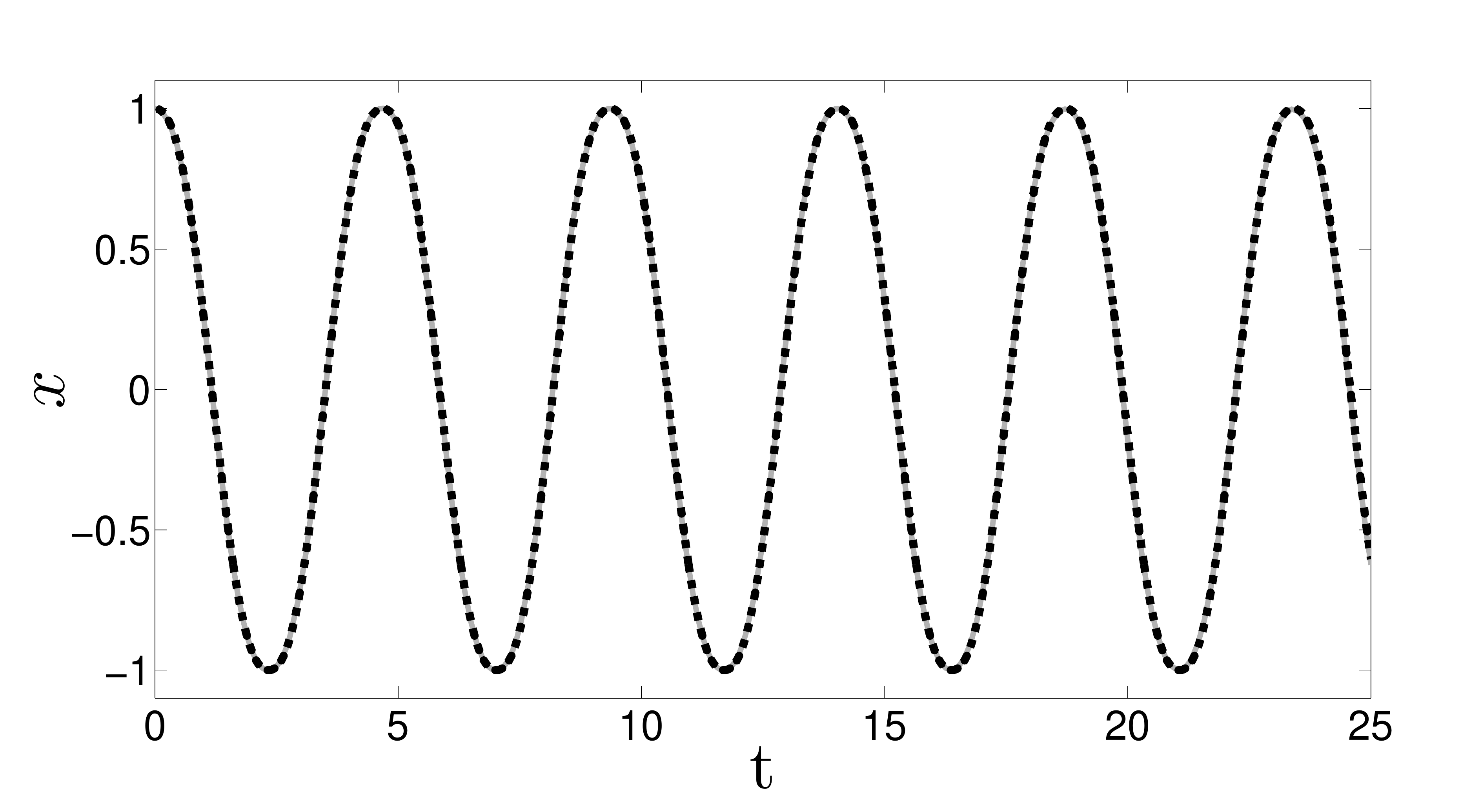}}
  \caption{Comparison of periodic solution obtained by HG and ode$45$ for $A=1,$ $\epsilon=3.2;$ continuous line: ode$45$, dashed line: HG.}
  \label{mod_HG}
\end{figure}

\begin{figure*}[h]
\centering{\includegraphics[scale=0.33]{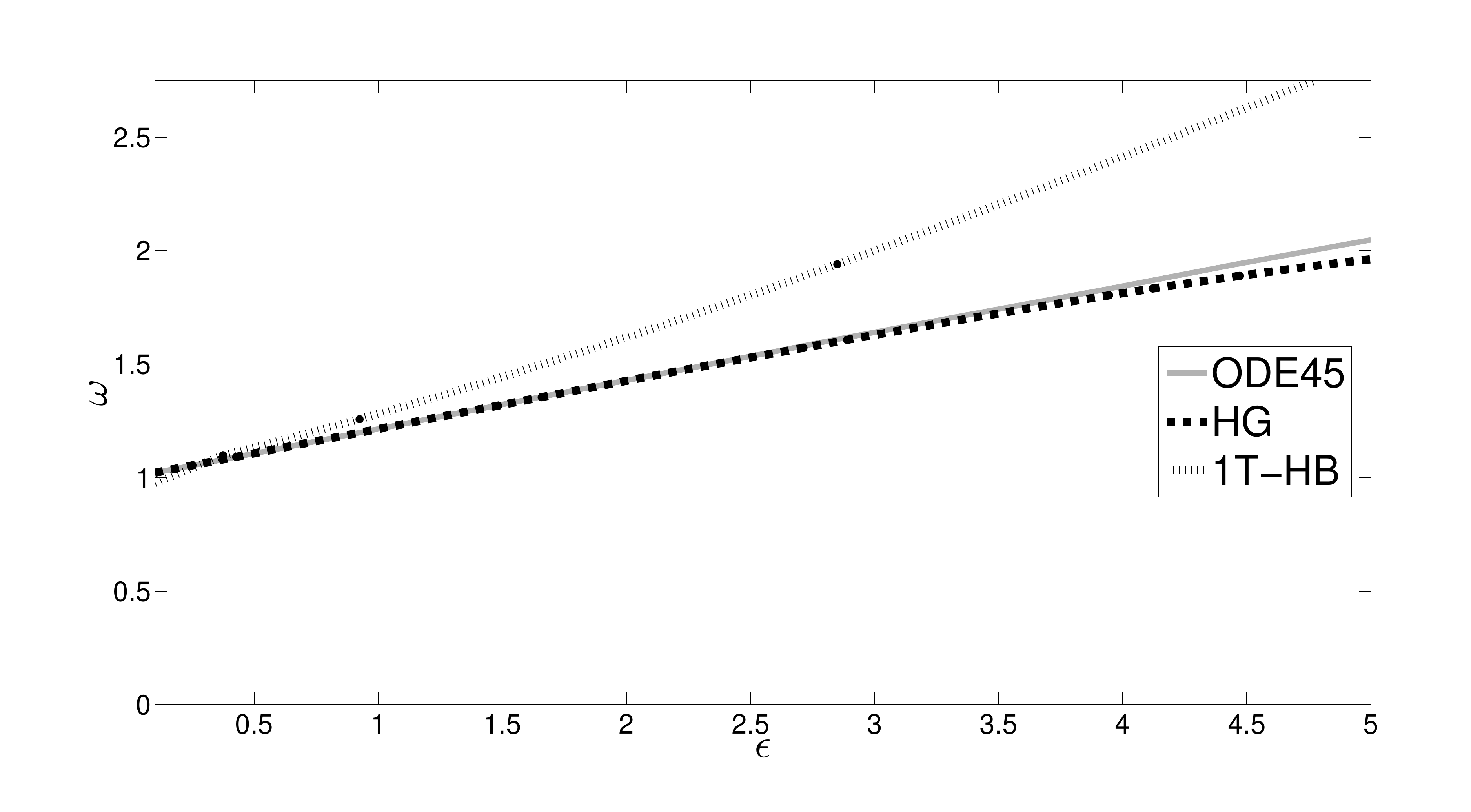}}
\caption{Natural frequency $\omega$ obtained for $A=1$ for oscillator, Eqn.(\ref{ODE_abs_simplified}).}
\label{freq_compare_abs}
\end{figure*}

We now consider an oscillator
\begin{equation}\label{ODE_abs}
\ddot{x} + x + \epsilon x |\dot{x}| = 0.
\end{equation}
The non-smoothness is now of type modulus with an argument of $\dot{x}.$ We re-write the above equation in terms of Heaviside function as
\begin{equation}\label{ODE_abs_simplified}
\ddot{x} + x + \epsilon x \dot{x}(2H(\dot{x})-1) = 0.
\end{equation}
The change in the argument of Heaviside from $x$ to $\dot{x}$ does not alter the framework discussed in Sec.(\ref{HG:heaviside}). Homotopy $\mathcal{H}$ in the scaled time $\tau$ is given by
\begin{equation*}\label{PDE_heaviside_scaled_abs}
\mathcal{H} \equiv (1-p) \Big(\tilde{x}_{\tau\tau} + \lambda(p)^2 \tilde{x}\Big) -h(p)\Bigg(\tilde{x}_{\tau\tau} + \lambda(p)^2 \tilde{x} + \epsilon\lambda(p)\tilde{x} \tilde{x}_{\tau}\left(2 H\left(\frac{\tilde{x}_{\tau}}{\lambda(p)}\right) - 1 \right)  \Bigg)=0.
\end{equation*}
Since the solutions to Eqn.(\ref{ODE_abs}) are not scalable, we obtain a periodic solution with initial condition $(A,0).$ We compute upto the first order and obtain approximate expressions for the natural frequency and solution $x_1(t)$ as 
\begin{dmath*}
\omega(1) = 1-\frac{\epsilon A}{3\pi}\left( {{4h^{\mbox{\tiny[1]}}}}+{{h^{\mbox{\tiny[2]}}}} +{{2{h^{\mbox{\tiny[1]}}}^{2}}}\right) -\epsilon^{2}A^2\left(\frac{{h^{\mbox{\tiny[1]}}}^{2}}{24}-\frac{\left( {{4h^{\mbox{\tiny[1]}}}}+{{h^{\mbox{\tiny[2]}}}} +{{2{h^{\mbox{\tiny[1]}}}^{2}}}\right)^2}{9\pi^2}\right),
\end{dmath*}
\begin{dmath*}
x_1(\tau) =  A\cos\tau + \epsilon A^2 h^{\mbox{\tiny[1]}}\left(\left( \frac{2}{3} \sin \tau+ \frac{1}{3} \sin 2\tau\right) H \left(\tau-\pi \right) + \left( \frac{2}{3} \sin \tau- \frac{1}{3} \sin 2\tau  \right)H \left(\tau-2\pi\right) + \frac{(2\pi-4\tau)}{6\pi}\sin \tau - \frac{1}{6}\sin\ 2\tau\right),
\end{dmath*}
where $\tau=\omega(1)t.$ We obtain the values of free parameters $h^{\mbox{\tiny[1]}}$ and $h^{\mbox{\tiny[2]}}$ by applying Galerkin projections with weighting functions $\cos(\omega(1)t)$ and $t \sin(\omega(1)t).$ We determine these unknowns for two values of $\epsilon$ (numerical values in Table \ref{table:unknowns:sig_abs}). We thus obtain $x_1(t)$ and compare it with the solution obtained via ode$45$ in Fig.(\ref{mod_HG}). Here as well, the expression for $x_1(t)$ has zero-crossings at $t=\frac{\pi}{2\omega(1)},\ {\frac{3\pi}{2\omega(1)}}$ making the evaluation of $H(x_1(t))$ easy.
Figure (\ref{freq_compare_abs}) shows the frequency obtained using different methods. Compared to the previous oscillator Eqn.(\ref{ODE_signum}), the $\epsilon$-usability range corresponding to the first order HG is restricted only to $\epsilon\in[0,4].$ The same in case of $1$T$-$HB is $\epsilon\in[0,1].$ It seems that the $\epsilon$-usability range is a strong function of the inherent nature of the oscillator, characterizing which is beyond the scope of current work.

\subsection{Inverse modulus type discontinuity}
\begin{figure}[!h]
	\centering
	{\includegraphics[draft=\status ,height=4cm,width=8cm]{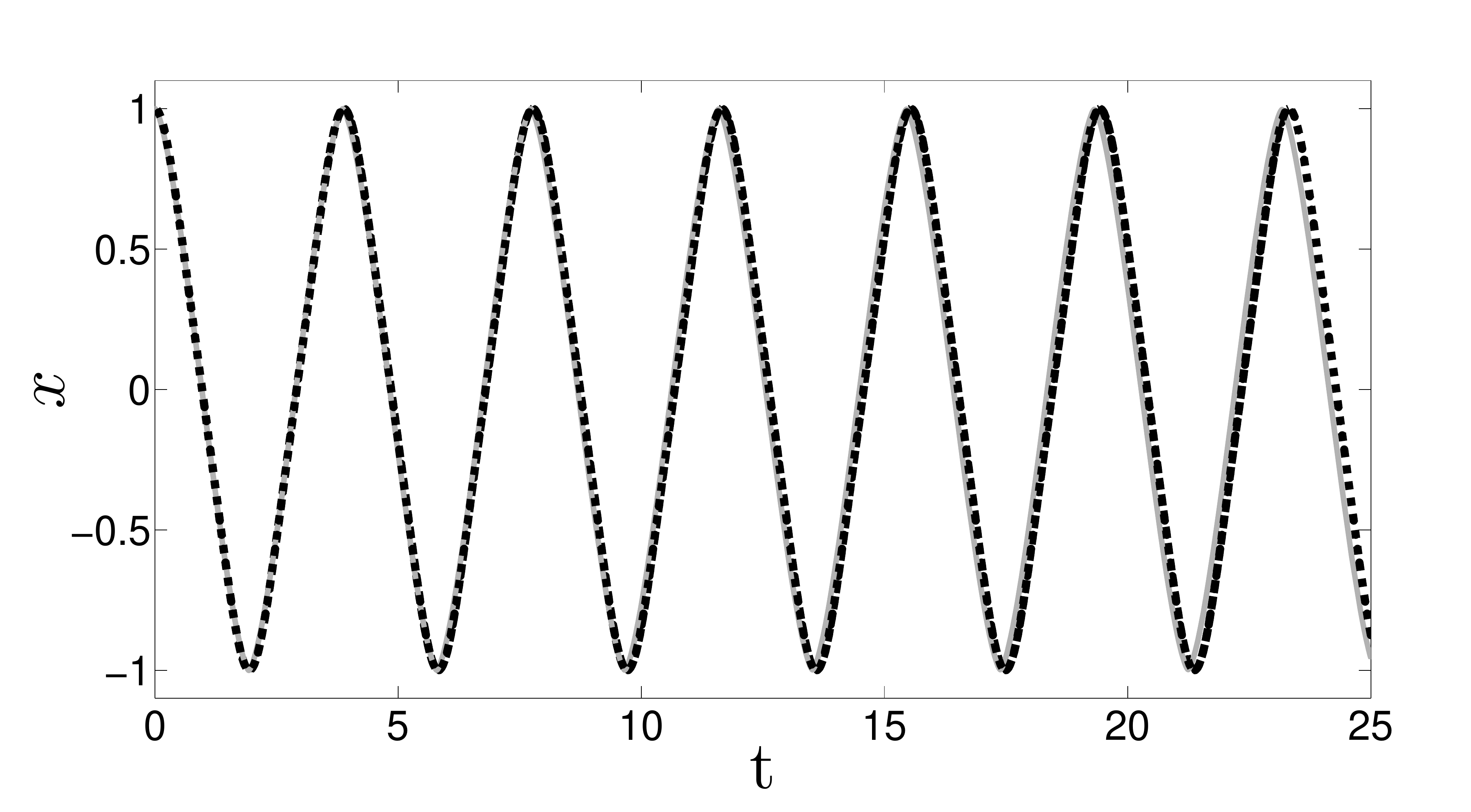}}
	\caption{Comparison of periodic solution obtained by HG and ode$45$ for $A=1,$ $\epsilon=1.6;$ continuous line: ode$45$, dashed line: HG.}
	\label{inv_mod_HG}
\end{figure}

Now, we consider an oscillator
\begin{equation}\label{ODE_abs_inv}
\ddot{x} + \epsilon\frac{x}{|\dot{x}|}=0.
\end{equation}
Multiplying the above equation by $|\dot{x}|$ and using $|\dot{x}|=\dot{x}(2H(\dot{x})-1),$ we get
\begin{equation}\label{ODE_abs_inv_simplified}
\dot{x}\ddot{x}(2H(\dot{x})-1) + \epsilon x = 0.
\end{equation}
A solution to this oscillator starting with an initial condition, say $(A,0)$ is periodic. 
The framework discussed in Sec.(\ref{HG:heaviside}) can be applied to Eqn.(\ref{ODE_abs_inv_simplified}) as well. The homotopy $\mathcal{H}$ in the scaled time $\tau$ is given by
\begin{equation*}
\mathcal{H} \equiv \lambda(p)(1-p) \Big(\tilde{x}_{\tau\tau} + \lambda(p)^2 \tilde{x}\Big) -h(p)\Bigg( \tilde{x}_{\tau} \tilde{x}_{\tau\tau} \left(2 H\left(\frac{\tilde{x}_{\tau}}{\lambda(p)}\right) - 1 \right) + \epsilon\lambda(p)^3 \tilde{x}  \Bigg)=0.
\end{equation*}
We obtain the expressions for the non-scalable periodic solution of Eqn.(\ref{ODE_abs_inv}) and corresponding frequency by applying HG upto the first order to get
\begin{dmath*}
\omega(1) = 1+\frac{A}{3\pi}\left({4{{h^{\mbox{\tiny[1]}}}}}+{{h^{\mbox{\tiny[2]}}}}\right)+\frac{A^2{h^{\mbox{\tiny[1]}}}^{2}}{12}- \epsilon\left( h^{\mbox{\tiny[1]}}+\frac{h^{\mbox{\tiny[2]}}}{4} \right) -\frac{3{\epsilon}^{2}	{h^{\mbox{\tiny[1]}}}^{2}}{8}
 \quad \mbox{and}
\end{dmath*}
\begin{dmath*}
x_1(\tau) =  A\cos\tau - A^2 h^{\mbox{\tiny[1]}}\left(\left( \frac{2}{3} \sin \tau+ \frac{1}{3} \sin 2\tau\right) H \left(\tau-\pi \right) + \left( \frac{2}{3} \sin \tau- \frac{1}{3} \sin 2\tau  \right)H \left(\tau-2\pi\right) + \frac{(\pi-2\tau)}{3\pi}\sin \tau - \frac{1}{6}\sin 2\tau\right),
\end{dmath*}
\[\mbox{with} \quad \tau=\omega(1)t. \]
The free parameters $h^{\mbox{\tiny[1]}}$ and $h^{\mbox{\tiny[2]}}$ are determined by applying the Galerkin projections, taking weighting functions as before, i.e., $\cos(\omega(1)t)$ and $t \sin(\omega(1)t).$ Computing the unknowns for two values of $\epsilon$ (numerical values, Table (\ref{table:unknowns:sig_abs})), in Fig.(\ref{inv_mod_HG}) we compare the solution so obtained with the one via MATLAB built-in integrator ode$45.$
\begin{table}[h]
\caption{Frequency and free parameters}
\renewcommand{\arraystretch}{1.25}%
\centering 
\begin{tabular}{|c|c| c c c| } 
\hline 
Equation & $\epsilon$ & $\omega(1)$ & $h^{\mbox{\tiny[1]}}$ & $h^{\mbox{\tiny[2]}}$ \\  
\hline  
{\multirow{2}{*}{(\ref{ODE_signum})}}  & $3.2$ &  $1.63719\cdots$ & $-0.36174\cdots$ & $-0.07434\cdots$\\ 
     & $32$ & $5.17727\cdots$ & $-0.03617\cdots$ &  $-0.03365\cdots$\\ [0.5ex] \hline
{\multirow{2}{*}{(\ref{ODE_abs})}} & $1.6$ &   $1.34203\cdots$ & $-0.85251\cdots$ &   $-0.00140\cdots$\\
     & $3.2$ & $1.66725\cdots$ & $-0.76351\cdots$ &  $-0.02312\cdots$\\ [0.5ex] \hline
{\multirow{2}{*}{(\ref{ODE_abs_inv})}} & $0.8$ &   $1.28342\cdots$ & $-1.73908\cdots$ &   $-0.44176\cdots$\\
     & $1.6$ & $1.61700\cdots$ & $-1.73908\cdots$ &  $-3.36354\cdots$\\[0.5ex] 
\hline 
\end{tabular}
\label{table:unknowns:sig_abs}
\end{table}                 

\begin{figure*}[h]
\centering{\includegraphics[scale=0.33]{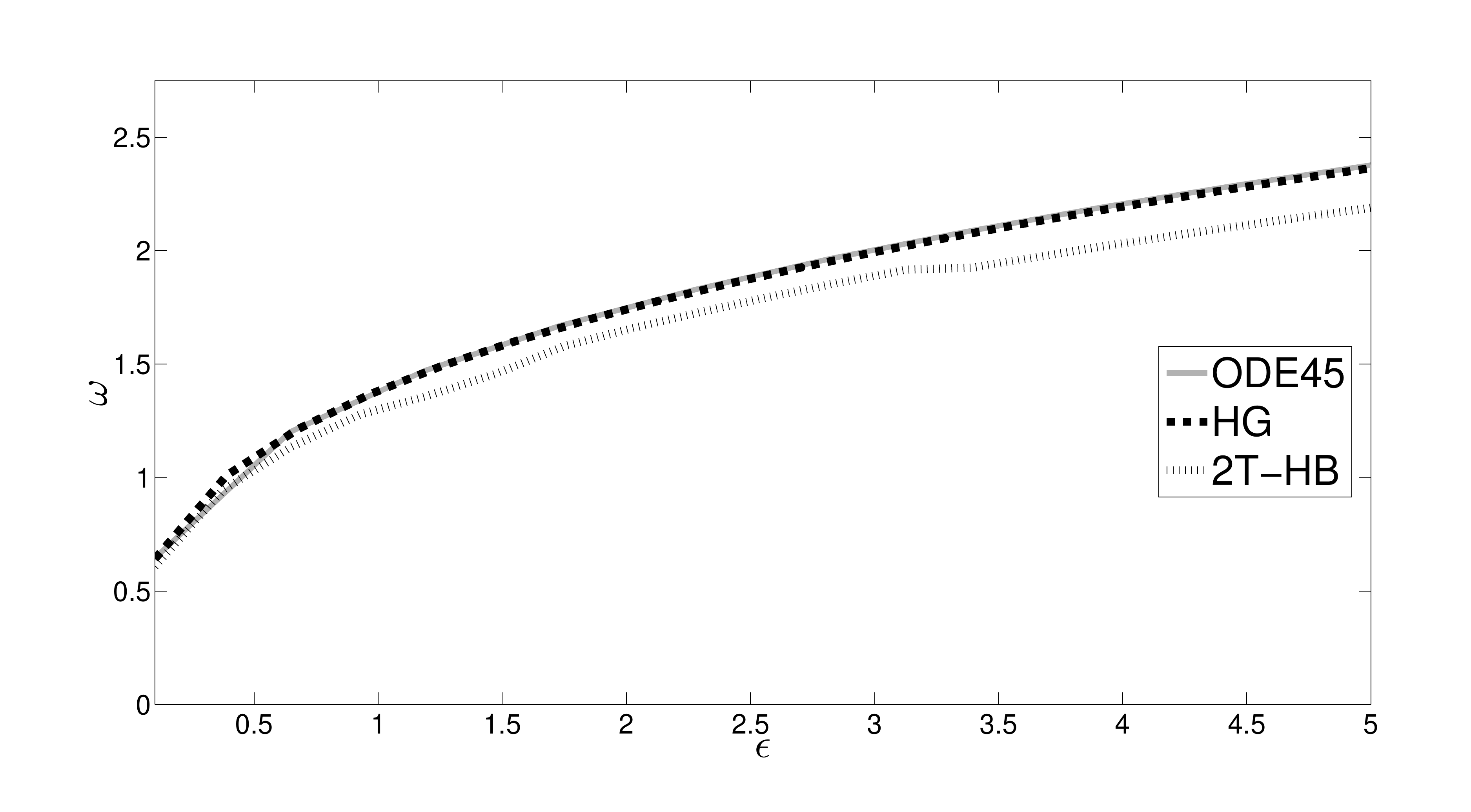}}
\caption{Natural frequency $\omega$ obtained for $A=1$ for oscillator, Eqn.(\ref{ODE_abs_inv_simplified}).}
\label{freq_compare_abs_inv}
\end{figure*}

Here again, the evaluation of $H(x_1(t))$ is easy as $x_1(t)$ has zero-crossings exactly at $t=\frac{\pi}{2\omega(1)},\ {\frac{3\pi}{2\omega(1)}}.$ Figure (\ref{freq_compare_abs_inv}) compares the frequencies by ode$45,$ HG and 2T-HB for $\epsilon\in[0.1,5].$ The $\epsilon$-usability range w.r.t. HG is $\epsilon=[0,10]$ and the same w.r.t. 2T-HB is $\epsilon=[0,0.8].$

With respect to all oscillators, HG method is consistently superior to other methods.

\section{Natural response of an impact oscillator}
\label{sec::5}
\subsection{The system description}

\begin{figure}[!h]
\centering{\includegraphics[scale=0.3]{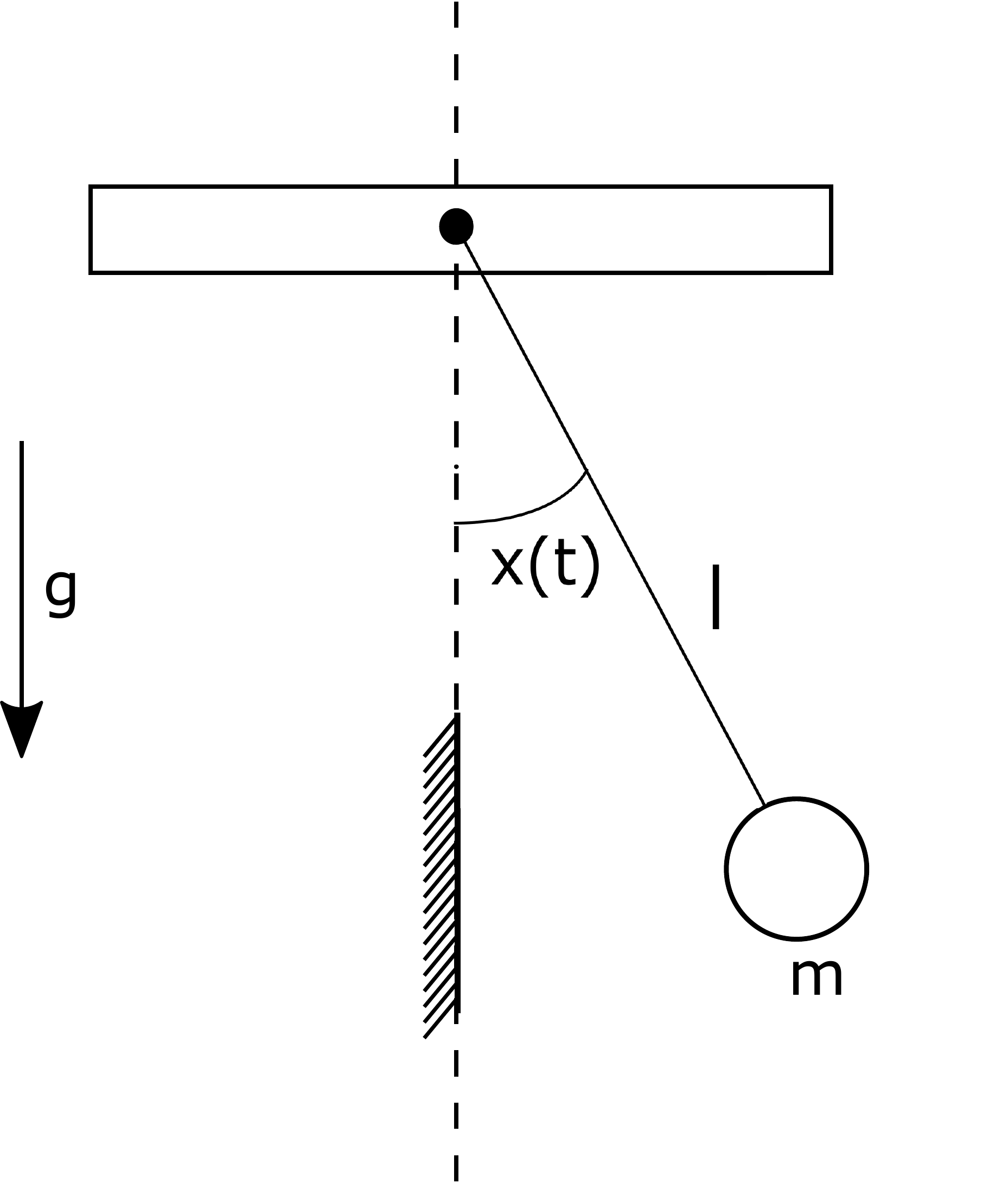}}
\caption{Pendulum with a vertical constraint}
\label{impact_pendulum}
\end{figure}

Consider a pendulum constrained unilaterally as shown in Fig.(\ref{impact_pendulum}). 
Scaling time by $\sqrt{\tfrac{g}{ml}}$ and using the Newtonian impact law, i.e.,
\[x(t^*)=0:\quad \dot{x}^+(t^*)= -e\dot{x}^-(t^*), \]
where $t^*$ is the instant of impact, $e$ is the coefficient of restitution, $\dot{x}^-(t^*)$ and $\dot{x}^+(t^*),$ denoting the velocity of the pendulum just before and after the collision respectively.

In order to apply homotopy analysis method, we wish to have an equation governing the motion of the pendulum along with the collision constraint of the form 
\begin{equation}\label{ODE_dirac_pre}
\ddot{x} + x + F_{\epsilon}(x, \dot{x}) = 0, 
\end{equation}
where $F_{\epsilon}(x, \dot{x})$ is expected to be a non-smooth function of $x.$
Introducing $\epsilon=-(1+e),$ we re-write the collision constraint as 
\begin{equation}\label{collision_constraint}
 \dot{x}^+(t^*) - \dot{x}^-(t^*) - {\epsilon}\dot{x}^-(t^*)=0.
\end{equation}
The displacement $x(t)$ increasingly away from the constraint is assumed to be positive. Then, velocity just before impact $\dot{x}^-(t^*)$ is negative. We thus modify the above taking $\dot{x}^-(t^*)=-|\dot{x}^-(t^*)|$ as 
\begin{equation*}
\dot{x}^+(t^*) - \dot{x}^-(t^*)+\epsilon\frac{\dot{x}^-(t^*)^2}{|\dot{x}^-(t^*)|}=0.
\end{equation*}         
Taking an infinitesimal positive number $\mu,$ $\dot{x}^-(t^*-\mu)$ and $\dot{x}^+(t^*+\mu),$ in the limit $\mu \rightarrow 0,$ represents the velocity of the pendulum just before and just after the impact respectively. Assuming $x(t)$ to be $\mathcal{C}^0$-continuous at the impact, we assume the left hand side of the above equation as
\begin{dmath}\label{ad-hoc}
 \lim_{\mu\to 0}\left(\dot{x}^+(t^*+\mu) - \dot{x}^-(t^*-\mu) + \int_{t^*-\mu}^{t^*+\mu}x\ dt+\epsilon \frac{\dot{x}^-(t^*)^2}{|\dot{x}^-(t^*)|} \right)=0.
\end{dmath}     
In order to formulate the governing law (\ref{ODE_dirac_pre}), we find it necessary to introduce a Dirac-delta like distribution, $\delta^*(t)$ defined as
\begin{equation*}
\int_{-\infty}^{\infty}x(t) \delta^*(t)dt=x^-(0),
\end{equation*}
where $x(t)$ is a non-smooth function and $x^-(0)$ is its left hand limit at $t=0.$ The properties of $\delta^*(t)$ are given in the appendix \ref{appendix1} with appropriate justification. Taking $y(t)=x(t)$ in the Eqn.(\ref{dirac_prop_new}) with the integral restricted to $[t^*-\mu,t^*+\mu]$ (with the time of impact $t^*$ being the root of $x$), Eqn.(\ref{ad-hoc}) becomes equivalent to 
\begin{dmath*}
 \lim_{\mu\to 0}\bigg(\dot{x}^+(t^*+\mu) - \dot{x}^-(t^*-\mu) + \int_{t^*-\mu}^{t^*+\mu}x\ dt+ \epsilon\int_{t^*-\mu}^{t^*+\mu} \dot{x}^2\delta^*(x)\ dt\bigg)=0.
\end{dmath*}
Re-writing the first two terms as the integral of $\ddot{x},$ we get
\begin{dmath*}
 \lim_{\mu\to 0}\bigg( \int_{t^*-\mu}^{t^*+\mu}\ddot{x}\ dt + \int_{t^*-\mu}^{t^*+\mu}x\ dt + \epsilon\int_{t^*-\mu}^{t^*+\mu} \dot{x}^2\delta^*(x)\ dt\bigg)=0
\end{dmath*}
 Thus,\[\lim_{\mu\to 0} \int_{t^*-\mu}^{t^*+\mu}\bigg(\ddot{x}\  + x+ \epsilon \dot{x}^2\delta^*(x)\bigg)\ dt=0. \]
The physics of the Newtonian impact suggests to hypothesize that the solution $x(t)$ is  $\mathcal{C}^0$-regular, which enables us to write the above  formally as 
\begin{equation}\label{ODE_dirac}
\ddot{x} + x + \epsilon \dot{x}^2 \delta^*(x) = 0 , 
\end{equation}
which when compared to Eqn.(\ref{ODE_dirac_pre}) yields $F_{\epsilon}(x, \dot{x})= \epsilon \dot{x}^2 \delta^*(x).$ The above equation is written formally and has to be understood in the sense of distributions, as precised in the appendix \ref{appendix1}. Equation (\ref{ODE_dirac}) holds for any sign convention for $x(t).$ 

Since the solutions to Eqn.(\ref{ODE_dirac}) are scalable, the frequency is  independent of initial conditions and is equal to that of the unconstrained version ($\epsilon=0$), i.e., one. This information plays an important role in application of HG as will be shown in the next subsection.

\subsection{Homotopy analysis and Galerkin projections}

Homotopy analysis method is typically used to find a periodic solution. Here we apply the method to obtain a decaying solution to Eqn.(\ref{ODE_dirac}) for some $-2\leq\epsilon< -1.$ Unlike an oscillator with viscous damping where the decay rate is exponential, here the same is algebraic. We construct the homotopy with linear operator as the viscously damped linear oscillator and proceed to formulate the homotopy as
\begin{align*}
    \mathcal{H} &\equiv (1-p)\mathcal{L} - h(p)\mathcal{N}=0, \\
\mbox{where}\quad	\mathcal{L} &\equiv \ddot{x} + \gamma(p) \dot{x} + \chi(p) x \quad \mbox{and} \quad \\
	 \mathcal{N} &\equiv \ddot{x} + x + \epsilon \dot{x}^2 \delta^*(x).
\end{align*}
$\gamma(p)$ and $\chi(p)$ are damping and stiffness as functions of $p,$ introduced to adjust the frequency and decay rate of the solution to $\mathcal{H}.$ Here, we strongly emphasize that the homotopy is written in the sense of distributions. Taking the assumed regularity of the solution to the homotopy at $p=1$ into account, we assume that $\tilde{x}(t;p)$ is a $\mathcal{C}^0$-regular function, which results in all terms of the homotopy including $\dot{\tilde{x}}(t;p)^2\delta^*(\tilde{x}(t;p))$ being distributions. Differentiation of each term w.r.t. $p$ is then justified in a distributional sense. We Taylor-expand the functions $\gamma(p)$ and $\chi(p)$ as 
\begin{align}
	\gamma(p) = \gamma(0) + \sum_{n=1}^{\infty}\frac{1}{n!} \gamma^{\mbox{\tiny [{\em n}]}}p^{n} \quad &\mbox{where} \quad \gamma^{\mbox{\tiny [{\em n}]}}=  \left.\frac{d^{n}\gamma(p)}{dp^{n}} \right|_{p=0},\label{gamma_expansion}\\
	\chi(p) =  \chi(0) + \sum_{n=1}^{\infty}\frac{1}{n!} \chi^{\mbox{\tiny [{\em n}]}}p^{n} \quad &\mbox{where} \quad \chi^{\mbox{\tiny [{\em n}]}} =  \left.\frac{d^{n}\chi(p)}{dp^{n}} \right|_{p=0}.\label{chi_expansion}
\end{align}

\begin{figure*}[!h]
\centering
{    \includegraphics[draft=\status ,height=6cm,width=12cm]{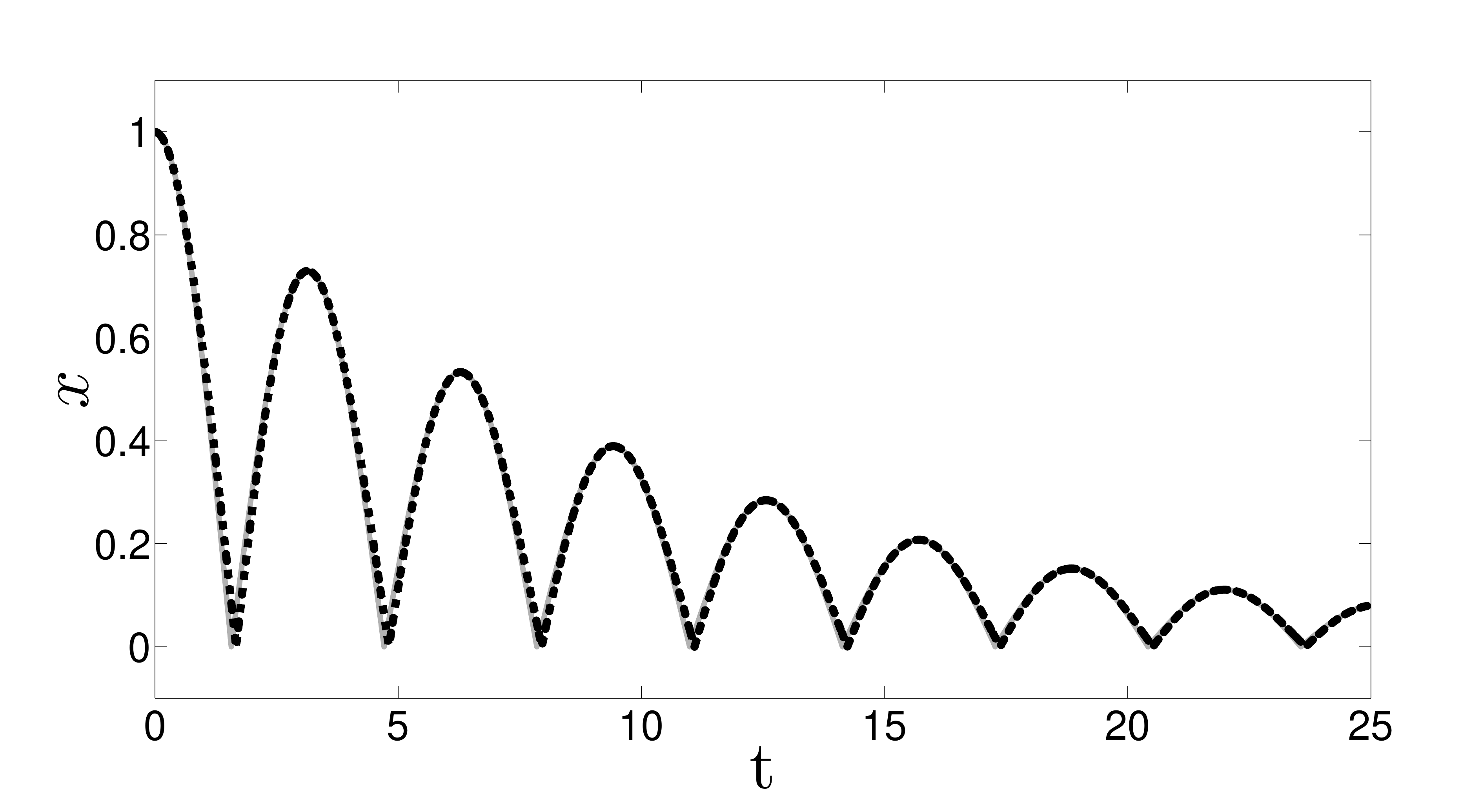}}

  \caption{Comparison of periodic solution obtained by HG and ode$45$ for $\epsilon=-1.73..$; continuous line: ode$45$, dashed line: HG.}
  \label{dirac_HG}
\end{figure*}

Considering the scalability of Eqn.(\ref{ODE_dirac}), we wish to obtain a solution with initial condition $(1,0).$ We begin by substituting $p=0$ in the homotopy and obtain the zeroth order deformation equation
\begin{equation}\label{zeroth_order_dirac}
\ddot{x}^{\mbox{\tiny [0]}} + \gamma(0) \dot{x}^{\mbox{\tiny [0]}} + \chi(0) x^{\mbox{\tiny [0]}} = 0.
\end{equation}
As can be seen from the above, ${\gamma(0)}$ is the decay rate and the expression $\omega_d(0)=\sqrt{\chi(0)-\frac{\gamma(0)^2}{4}}$ is the damped natural frequency of the oscillator at the zeroth order. However, this true for all other orders as well. Hence ${\gamma(0)}$ is the decay rate and $\omega_d(0)$ is the frequency of the solution to the homotopy at $p=1,$ the desired solution. $\gamma(0)$ and $\chi(0)$ can be determined by Galerkin projections; however, this not only increases the number of equations to be solved, but also makes the equations non-algebraic in the unknowns ($\gamma(0),\ \chi(0),\ h^{\mbox{\tiny [1]}},\ h^{\mbox{\tiny [2]}},\cdots$ ) thereby complicating the application of Galerkin projections. Hence, we prefer to pre-determine $\gamma(0)$ and $\chi(0)$ based on the physics of the impact oscillator as follows. As discussed previously, the frequency of the solution to the unilaterally constrained pendulum is equal to that of the unconstrained one, therefore we set 
\begin{equation*}
\omega_d(0)=\sqrt{\chi(0)-\frac{\gamma(0)^2}{4}}=1.
\end{equation*}
Equating the decay rate per cycle of the zeroth order deformation equation to that of Eqn.(\ref{ODE_dirac}) and then from the above expression, we obtain
\begin{equation}\label{gamma0_chi0}
\gamma(0) ={\frac {-2\ln  \left( -\epsilon-1 \right) }{\sqrt {  [\ln  \left( -\epsilon-1 \right)]^{2}+4 {\pi }^{2}}}} \quad \mbox{and}\quad 
  \chi(0) =\sqrt{1+{\frac {[\ln  \left( -\epsilon-1 \right)]^2 }{ [\ln  \left( -\epsilon-1 \right)]^{2}+4 {\pi }^{2}}}}.
\end{equation}
Choosing the initial condition as $(1,0)$ for Eqn.(\ref{zeroth_order_dirac}), we solve it and obtain
\begin{equation}\label{first_order_dirac}
x^{\mbox{\tiny [0]}}(t) = \exp(-\tfrac{\gamma(0) t}{2}) \left(\cos t + \frac{\gamma(0)}{2}\sin t\right),
\end{equation}
where $\gamma(0)$ is given by Eqn.(\ref{gamma0_chi0}). Zeroes of $x^{\mbox{\tiny [0]}}(t),$ important in the subsequent analysis are given by 
\begin{equation}\label{roots_x0}
t_i=i\pi-\tan^{-1}\left(\frac{2}{\gamma(0)}\right),\quad \mbox{where}\quad i=1,2,3,\cdots.
\end{equation}
Differentiating the homotopy once, substituting $p=0,$ and using Eqns.(\ref{taylor_x_LPP}),(\ref{solution_taylor_h}),(\ref{gamma_expansion}),(\ref{chi_expansion}), we obtain the first order deformation equation,
\begin{equation*}
\ddot{x}^{\mbox{\tiny [1]}} + \gamma(0) \dot{x}^{\mbox{\tiny [1]}} + \chi(0) x^{\mbox{\tiny [1]}} = h^{\mbox{\tiny [1]}}\Big(\ddot{x}^{\mbox{\tiny [0]}} + x^{\mbox{\tiny [0]}} + \epsilon {x^{\mbox{\tiny [0]}}}^2 \delta^*(x^{\mbox{\tiny [0]}}) \Big) - \gamma^{\mbox{\tiny [1]}} \dot{x}^{\mbox{\tiny [0]}} - \chi^{\mbox{\tiny [1]}} x^{\mbox{\tiny [0]}}=:F_1(t).
\end{equation*}
where $x^{\mbox{\tiny [0]}}$ is given by Eqn.(\ref{first_order_dirac}). To obtain a solution to the above with ease, we simplify the term $\delta^*(x^{\mbox{\tiny [0]}}(t)).$ For any function $f(t)$ of class $\mathcal{C}^1,$ Dirac-delta satisfies
\begin{equation}\label{dirac_prop}
 \delta(f(t))= \sum_{i=-\infty}^{\infty}\frac{\delta(t-t_i)}{|\dot{f}(t_i)|},
 \end{equation}
 where $t_i$ are the roots of $f(t)$ for $t>0.$
 Since $x^{\mbox{\tiny [0]}}(t)$ is smooth, Dirac-delta distribution property \eqref{dirac_prop} is satisfied by $\delta^*(x^{\mbox{\tiny [0]}}(t)).$ Using this equality and taking $f(t)= x^{\mbox{\tiny [0]}}(t)$ in Eqn.(\ref{dirac_prop}),  we get
\begin{equation}\label{dirac_star_expansion}
\delta^*(x^{\mbox{\tiny [0]}}(t)) =  \sum_{i=-\infty}^{\infty}\frac{\delta(t-t_i)}{\sqrt{(1+\frac{\gamma(0)^2}{4})}\exp(-\tfrac{\gamma(0)t_i}{2})}, 
\end{equation}
where $t_i$'s are given by Eqn.(\ref{roots_x0}). Substituting the above expression in the first order deformation equation, we get an oscillator forced at its damped natural frequency, $1$. We remove terms from the forcing which have the frequency equal to that of the unforced oscillator, even though amplitudes of such terms are not constant w.r.t. time. Therefore, equating the coefficients of $\sin t$ and $\cos t$ in $F_1(t)$ to zero, we get two equations algebraic in unknowns $\gamma^{\mbox{\tiny [1]}}$ and $\chi^{\mbox{\tiny [1]}},$ which on solving yield
\[ \gamma^{\mbox{\tiny [1]}}= -\frac{\gamma(0)^2 h^{\mbox{\tiny [1]}}}{\gamma(0)^2+{4}}\quad \mbox{and} \quad \chi^{\mbox{\tiny [1]}} = -\gamma(0)h^{\mbox{\tiny [1]}}.\]
These expressions play no role in further analysis.
The first order deformation equation thus becomes
\begin{equation*}
\ddot{x}^{\mbox{\tiny [1]}} + \gamma(0) \dot{x}^{\mbox{\tiny [1]}} + \chi(0) x^{\mbox{\tiny [1]}} = h^{\mbox{\tiny [1]}}\epsilon {x^{\mbox{\tiny [0]}}}^2 \delta^*(x^{\mbox{\tiny [0]}})=:F_1(t),
\end{equation*}
where $\delta^*(x^{\mbox{\tiny [0]}})$ is as expressed in Eqn.(\ref{dirac_star_expansion}). The above can be re-written as a system of two first order differential equations such that the homogeneous part  satisfies the Carathéodory conditions and the non-homogeneous part is a distribution that is not a Lebesgue integrable, function. As per the analysis from chapter 1 of \cite{Filipov}, it follows that the first order equation admits a unique solution of class $\mathcal{C}^0.$ We then solve it with initial condition $(0,0)$ and obtain
\[x^{\mbox{\tiny [1]}}(t)= \epsilon h^{\mbox{\tiny [1]}} x^{\mbox{\tiny [0]}}(t)\sum_{i=1}^{\infty}(-1)^i H\left(t-t_i\right), \]
where $H$ denotes the Heaviside function and $t_i$'s are given by Eqn.(\ref{roots_x0}). Adding $x^{\mbox{\tiny [0]}}(t)$ and $x^{\mbox{\tiny [1]}}(t),$ we get the expression for the solution $x_1(t).$ We now determine the only unknown $h^{\mbox{\tiny [1]}}$ by applying Galerkin projections with weighting function $w_1(t)=\cos t\sum_{i=1}^{\infty}H\left(t-t_i\right).$ Computationally, it suffices to consider the first four zeroes of $x^{\mbox{\tiny [0]}}(t)$ in $w_1(t)$ to get a good match with the solution obtained via numerical integration. We obtain the residual by substituting $x_1(t)$ in Eqn.(\ref{ODE_dirac}) as
\begin{equation}\label{residual_dirac}
\mathcal{R}_1(t)= \ddot{x}_1 + x_1 + \epsilon \dot{x}_1^2 \delta^*(x_1),
\end{equation}
and minimize the weighted residual for all time $t$ as
\begin{equation}\label{galerkin_dirac}
	\int_{-\infty}^{\infty} w_1(t)\mathcal{R}_1(t) dt = 0.
\end{equation}
In order to ease the integration of the above, we use the Eqn.(\ref{dirac_prop_new}) with $y(t)=x_1(t)$ and $f(y(t),\dot{y}(t))=1$  and obtain 
\[ \delta^*(x_1(t))= \sum_{i=-\infty}^{\infty}\frac{\delta(t-t_i)}{|\dot{x}_1^-(t_i)|}.\]
Using the above, we evaluate the integral in Eqn.(\ref{galerkin_dirac}) and determine $h^{\mbox{\tiny [1]}}.$ 
Simpler choice of $\gamma(0)$ is nicer to implement in Maple than a simpler choice of $\epsilon.$ For $\gamma(0)=0.2,$ inverting Eqn.(\ref{gamma0_chi0}), we arrive at $\epsilon=-1.73...$.
In order to compute the numerical solution to Eqn.(\ref{ODE_dirac}), we use MATLAB's ode$45$ integrator with ``event detection" on. While integrating the simple harmonic oscillator, we detect the zero-crossings of $x(t).$ Every time the event is detected, we change the velocity $\dot{x}(t)$ according to Newtonian impact law and resume numerical integration of the oscillator with the changed velocity. For $\epsilon=-1.73...,$ the solution $x_1(t)$ is compared to the solution obtained using ode$45$ in Fig.(\ref{dirac_HG}). The choice of linear operator ensures that the decay rate as well as the frequency are captured well via HG method. But this also results in the shifting of the zero crossings of the HG solution by some fixed amount from that of the numerical solution. We believe that this limitation maybe overcome by better choosing the linear operator. Although the global behaviour of the solution to Eqn.(\ref{ODE_dirac}) is like a linearly damped solution, its behaviour between two consecutive impacts is like an undamped solution; which is unlike the solution resulting from our choice of linear operator. An intelligent choice of the linear operator, which is partly though significantly, guided by the physics of the impact oscillator enables to get a good match with first two terms of the series solution. Practical application of HAM and thus HG demand the knowledge of the physical laws governing the problem.

\section{Conclusion and further work}
We have applied homotopy analysis method in combination with the Galerkin projections to several non-smooth oscillators and have obtained periodic responses approximately. The methodology developed for an oscillator with discontinuity of type Heaviside may easily be extended to discontinuities of type signum and modulus, as illustrated. Introduction of the convergence-control function and the time-stretching function make the framework natural and simpler. This renders us with an approximate analytical expression for the natural frequency involving the unknowns, also to be solved simultaneously but numerically. The natural frequency thus obtained matches excellently with the frequency obtained via numerical integration for a large $\epsilon$-range. The comparison establishes the superiority of homotopy analysis combined with Galerkin projections over purely heuristic method such as harmonic balance and purely asymptotic Lindstedt-Poincaré and NSTT+KB. With certain limitations, homotopy analysis has also been applied to the unilaterally constrained simple pendulum where the solutions are scalable, but decaying with time.
\par
With suitable modifications, we believe that the methodology developed here can be successfully applied to capture non-periodic response of a non-smooth oscillator, forced or unforced. We also believe that HG is useful not only in obtaining the natural frequencies of multi-DOF non-smooth oscillator systems but also in approximating solutions of boundary value problems with non-smooth governing equations. \\
\vspace{0.2in}

\textbf{Conflict of Interest} The authors declare that they have no conflict
of interest. \\
\vspace{0.2in}

\textbf{Acknowledgements} Heavy computational tasks involved in this work were supported through SERB-DST grant EMR/2014/001246.

\appendixpage
\begin{appendices}
\raggedbottom
\section{Properties of $\delta^*$ distribution}\label{appendix1}
Let $\Omega$ be an open set in $\mathbb{R}$ encompassing $0.$ By  $\mathcal{C}^{\infty}_c(\Omega),$ we denote the space of $\mathcal{C}^{\infty}$ regular functions with compact support in $\Omega.$ The members of the dual space of $\mathcal{C}^{\infty}_c(\Omega)$ are distributions. We say that an equation $f(x,\dot{x},\ddot{x})=0$ is satisfied in the sense of distribution if $$ \int_{\Omega}f(x,\dot{x},\ddot{x})\varphi\ dt=0\quad \forall \varphi\in \mathcal{C}^{\infty}_c(\Omega).$$

We introduce a Dirac-delta-like function $\delta^*(t)$ and the distribution generated by it $T_{\delta^*}(\cdot)$ as
\begin{equation}\label{dirac_prop_1}
T_{\delta^*}(\varphi)=\int_{\Omega}\varphi(t) \delta^*(t)dt=\varphi^-(0)\quad  \forall \varphi \in \mathcal{C}^{\infty}_c(\Omega),
\end{equation}
where $\varphi^-(0)$ is the left hand limit of the function $\varphi(t)$ at $t=0.$
In order to determine if $\delta^*(t)$ is a distribution, we take any convergent sequence $\{\varphi_n\}_{n\geq 1}\in \mathcal{C}^{\infty}_c(\Omega)$ converging to $\varphi\in \mathcal{C}^{\infty}_c(\Omega),$ and check $|T_{\delta^*}(\varphi_n)-T_{\delta^*}(\varphi)|=|\varphi_n^-(0)-\varphi^-(0)|=|\varphi_n(0)-\varphi(0)|\longrightarrow 0$ as $n\longrightarrow \infty.$ This implies that $T_{\delta^*}(\varphi_n)\longrightarrow T_{\delta^*}(\varphi) $ as $n\longrightarrow \infty$ and thus $\delta^*(t)$ is a distribution.

Let $\alpha$ be a whole number. The derivative $\partial^{\alpha} \delta^*(t)$ is given by 
\begin{equation}\label{dirac_prop_2}
T_{\partial^{\alpha}\delta^*}(\varphi)=\int_{\Omega}\varphi(t) \partial^{\alpha}\delta^*(t)dt=(-1)^{\alpha}\int_{\Omega}\partial^{\alpha} \varphi(t)\delta^*(t)dt=(-1)^{\alpha}\partial^{\alpha}\varphi^-(0)\quad  \forall \varphi \in \mathcal{C}^{\infty}_c(\Omega),
\end{equation}
where $\partial^{\alpha} \varphi^-(0)$ denotes the left hand limit of the function $ \partial^{\alpha} \varphi(t)$ at $t=0.$ The above property can be extended to $\varphi(t)\in \mathcal{C}^{\alpha}_c(\Omega),$ (space of $\alpha$-regular functions with compact support in $\Omega$) by a density argument. We remark that this property holds for the Dirac-delta distribution as well.

We observe that like the Dirac-delta distribution,  $\delta^*(t)$ satisfies the property 
\begin{equation}\label{dirac_prop_3}
\int_{-\infty}^{\infty} \delta^*(at)\ dt =\frac{1}{|a|} \int_{-\infty}^{\infty} \delta^*(t)\ dt \quad \forall a\in\mathbb{R}\setminus \{0\},
\end{equation}
by applying a change of variable $at\mapsto t.$

Let us consider a series $\{f_n\}_{n\geq 1}$ given by 
\[f_{n}(t)= \begin{cases} 
      \frac{2n}{\sqrt{\pi}}e^{-(nt)^2} & t\leq 0 \\
      0 & 0< t ,
   \end{cases}
\]
and check if $f_n(t)\longrightarrow \delta^*(t)$ as $n\longrightarrow \infty $ in the distributional sense. By $T_{f_n},$ we denote the distribution generated by $f_n$ on $\mathcal{C}^{\infty}_c(\Omega).$ For all $\varphi(t)\in \mathcal{C}^{\infty}_c(\Omega)$ extended by $0$ outside $\Omega,$ we have
\begin{equation}\label{dirac_prop_pre}
    T_{f_n}(\varphi) =\int_{-\infty}^{\infty}f_n(t)\varphi(t) dt=\int_{-\infty}^{0}\frac{2n}{\sqrt{\pi}}e^{-(nt)^2}\varphi(t)dt=\int_{-\infty}^{0}\frac{2}{\sqrt{\pi}}e^{-x^2}\varphi(x/n)dx,
\end{equation}
where we have applied a variable transform $t=x/n$ in the second last equality. The function series $\frac{2}{\sqrt{\pi}}e^{-x^2}\varphi(x/n)$ is bounded by an integrable function $\frac{2}{\sqrt{\pi}}e^{-x^2}\norm{\varphi}_{L^{\infty}(\Omega)}.$ Moreover, the series converges pointwise to $\frac{2}{\sqrt{\pi}}e^{-x^2}\varphi^-(0)$ for $x\in (-\infty,0]$ as $n\longrightarrow \infty.$ Hence applying dominated convergence theorem to the last integral in Eqn.(\ref{dirac_prop_pre}) gives
\begin{dmath}
\lim_{n\rightarrow\infty}T_{f_n}(\varphi) =\lim_{n\rightarrow\infty}\int_{-\infty}^{0}\frac{2}{\sqrt{\pi}}e^{-x^2}\varphi(x/n)dx = \int_{-\infty}^{0}\frac{2}{\sqrt{\pi}}e^{-x^2}\varphi^-(0)dx= \varphi^-(0).
\end{dmath}
Using the approximating function series $\{f_n\}_{n\geq 1}$, we calculate the integral of $\delta^*(t)$ in $t\in[a,0],\ a< 0,$ and find it exactly equal to its integral in $t\in[a,b],\ a<0<b.$ Or, equivalently, 
\begin{equation}\label{dirac_prop_4}
    \int_{-\infty}^{0}\varphi(t)\delta^*(t)\ dt=\int_{-\infty}^{\infty}\varphi(t)\delta^*(t)\ dt,\quad \forall \varphi\in \mathcal{C}^{\infty}_c(\Omega). 
\end{equation}

Now, we consider the distribution $\delta^*(x(t)),$ where $x(t)$ is a non-smooth function. Let time instant $t^*$ be such that $x(t^*)=0.$ We take an interval $I\equiv [t^*-e,t^*+e]$ (for some $e> 0$) such that for all $t\in I\setminus\{ 0\}, \ x(t)\neq 0.$ Assuming $\eta\in[0,e]$, for the left half of the interval $I,$ we have $$ x(t^*-\eta)=x(t^*)-\dot{x}^-(t^*)\eta +\mathcal{O}(\eta^2)=-\dot{x}^-(t^*)\eta +\mathcal{O}(\eta^2),$$ using which, we compute
\begin{dmath*}\label{final_prop}
\int_{t^*-e}^{t^*+e}\delta^*(x(t))\ dt = \int_{t^*-e}^{t^*}\delta^*(x(t))\ dt \hfill \mbox{(by Eqn.(\ref{dirac_prop_4}))}
= -\int_{e}^{0}\delta^*(x(t^*-\eta))\ d\eta
=-\int_{e}^{0}\delta^*(-\dot{x}^-(t^*)\eta +\mathcal{O}(\eta^2))\ d\eta
=\int_{-e}^{0}\delta^*(\dot{x}^-(t^*)\eta +\mathcal{O}(\eta^2))\ d\eta
=\int_{-e/|\dot{x}^-(t^*)|}^{0} \frac{\delta^*(\eta)}{|\dot{x}^-(t^*)|}\ d\eta \hfill\mbox{(by Eqn.(\ref{dirac_prop_3}))} = \frac{1}{|\dot{x}^-(t^*)|}=\int_{t^*-e}^{t^*+e}\frac{1}{|\dot{x}^-(t^*)|}\delta^*(t-t^*)\ dt.
\end{dmath*}
Generalizing the above property and using the definition of $\delta^*(t),$ we find that for some $f(y(t),\dot{y}(t))$ with non-smooth $y(t),$ we have
\begin{equation}\label{dirac_prop_new}
\int_{-\infty}^{\infty} f(y(t), \dot{y}(t))\delta^* (y(t)) dt = \sum_{i=-\infty}^{\infty}\frac{f^-(y(t_i), \dot{y}(t_i))}{|\dot{y}^-(t_i)|},
\end{equation}
where $t_i$'s are the roots of $y(t).$ Here as well, we remark that the above property holds in the case of Dirac-delta function when $y(t)$ is $\mathcal{C}^1$ regular and $f(\cdot,\cdot)$ is $\mathcal{C}^0$ regular.

\end{appendices}


\end{document}